\documentclass[journal]{IEEEtran}
\ifCLASSINFOpdf
\else
\fi
\usepackage{amssymb}
\usepackage{amsmath}

\usepackage{bbding}
\usepackage{amsfonts}
\usepackage{cite}
\usepackage{cases}
\usepackage{subfigure}
\usepackage{algorithmic}
\usepackage{mathrsfs}
\usepackage{bm}
\usepackage{verbatim}
\usepackage{enumerate}
\usepackage{booktabs}
\usepackage{multirow}
\usepackage[linesnumbered,ruled,longend]{algorithm2e}
\usepackage{textcomp}
\usepackage{slashbox}
\usepackage{graphics}
\usepackage{graphicx}
\usepackage{subeqnarray}
\usepackage{cases}

\usepackage{url}

\usepackage{array}
\newcommand{\PreserveBackslash}[1]{\let\temp=\\#1\let\\=\temp}
\newcolumntype{C}[1]{>{\PreserveBackslash\centering}p{#1}}
\newcolumntype{R}[1]{>{\PreserveBackslash\raggedleft}p{#1}}
\newcolumntype{L}[1]{>{\PreserveBackslash\raggedright}p{#1}}

\def\ba{\begin{array}}
\def\ea{\end{array}}
\newcommand{\beq}{\begin{equation}}
\newcommand{\eeq}{\end{equation}}
\newcommand{\bq}{\begin{eqnarray}}
\newcommand{\eq}{\end{eqnarray}}
\newcommand{\bqn}{\begin{eqnarray*}}
\newcommand{\eqn}{\end{eqnarray*}}
\newcommand{\bee}{\begin{enumerate}}
\newcommand{\eee}{\end{enumerate}}
\newcommand{\bi}{\begin{itemize}}
\newcommand{\ei}{\end{itemize}}

\IEEEoverridecommandlockouts

\usepackage{colortbl}

\usepackage{comment}
\newboolean{showcomments}
\setboolean{showcomments}{true}
\newcommand{\slow}[1]{\ifthenelse{\boolean{showcomments}}
{ \textcolor{red}{(SL:  #1)}}{}}
\newcommand{\you}[1]{\ifthenelse{\boolean{showcomments}}
{ \textcolor{green}{(PCY:  #1)}}{}}

\begin{document}

\title{Scheduling of EV Battery Swapping, I:  Centralized Solution}

\author{Pengcheng~You,~\IEEEmembership{Student~Member,~IEEE,}
        Steven~H.~Low,~\IEEEmembership{Fellow,~IEEE,}\\
        Wayes~Tushar,~\IEEEmembership{Member,~IEEE,}
        Guangchao~Geng,~\IEEEmembership{Member,~IEEE,}\\
        Chau~Yuen,~\IEEEmembership{Senior~Member,~IEEE,}
        Zaiyue~Yang,~\IEEEmembership{Member,~IEEE,}
        and Youxian~Sun
\thanks{P. You, G. Geng, Z. Yang and Y. Sun are with the State Key Laboratory of Industrial Control Technology, Zhejiang University, Hangzhou, 310027, China (e-mail: pcyou@zju.edu.cn; ggc@zju.edu.cn; yangzy@zju.edu.cn; yxsun@iipc.zju.edu.cn).}
\thanks{P. You and S. H. Low are with the Engineering and Applied Science Division, California Institute of Technology, Pasadena, CA 91125 USA (e-mail: pcyou@caltech.edu; slow@caltech.edu).}
\thanks{W. Tushar and C. Yuen are with Singapore University of Technology and Design (SUTD), Singapore 487372 (e-mail: wayes\_tushar@sutd.edu.sg; yuenchau@sutd.edu.sg).}
}

\maketitle

\begin{abstract}
We formulate an optimal scheduling problem for battery swapping that
assigns to each electric vehicle (EV) a best station to swap its depleted
battery based on its current location and state of charge.   The schedule aims
to minimize total travel distance and generation cost over both station
assignments and power flow variables, subject to EV range constraints,
grid operational constraints and AC power flow equations.   To deal with
the nonconvexity of power flow equations and the binary nature of station assignments, we propose a solution based on second-order cone programming (SOCP)
relaxation of optimal power flow (OPF) and generalized Benders decomposition.  When
the SOCP relaxation is exact, this approach computes a
globally optimal solution.  We evaluate the performance of the proposed algorithm
through simulations.  The algorithm requires global information and
is suitable for cases where the distribution network, stations,
and EVs are managed centrally by the same operator.  In Part II
of the paper, we develop distributed solutions for cases where they
are operated by different organizations that do not share private information.
\end{abstract}


\IEEEpeerreviewmaketitle


\section{Introduction}\label{sec:introduction}

\subsection{Motivation}

We are at the cusp of a historic transformation of our energy system into a more sustainable form in the coming decades. Electrification of our transportation system will be an important component because vehicles today consume more than a quarter of  energy in the US and emit more than a quarter of energy-related carbon dioxide \cite{C2ES2016climate, eia2015annual}. Electrification will not only greatly reduce greenhouse gas emission, but will also have a big impact on the future grid because electric vehicles are large but
flexible loads \cite{leou2014stochastic}.
It is widely believed that uncontrolled EV charging may stress the distribution
grid and cause voltage instability, but well controlled charging can help stabilize the
grid and integrate renewables.   As we will see below there is a large literature on
various aspects of EV charging.

We study a different problem here, motivated by a battery swapping
model currently being pursued in China, especially for electric buses and electric taxis
\cite{shang2015orchestrating}.
The State Grid (one of the two national utility companies) of
China is experimenting with a new business model where it operates not only the grid,
but also battery stations and a taxi service around a city, e.g., Hangzhou.  When the
state of charge of a State Grid taxi is low, it goes to one of State Grid operated battery
stations to exchange its depleted battery for a fully-charged battery.  While battery
swapping takes only a few minutes, it is not uncommon for a taxi to arrive at a station
only to find that it runs out of fully-charged batteries and there is a queue of taxis
waiting to swap their batteries.   The occasional multi-hour waits are a serious impediment
to the battery swapping model.

In this paper, we formulate in Section \ref{sec:system model}
an optimal scheduling problem for battery swapping that
assigns to each EV a best station to swap its depleted
battery based on its current location and state of charge.  The station assignment
not only determines EVs' travel distance, but can also impact significantly the
power flows on a distribution network because batteries are large loads.
The schedule aims to minimize a weighted sum of total travel distance and generation cost over both station assignments and power flow variables,
subject to EV range constraints, grid operational constraints and AC power flow equations.

This joint battery swapping scheduling and OPF problem
is nonconvex and computationally difficult for two reasons.
First the AC power flow equations are nonlinear.
Second, the station assignment variables are binary.
We address the first difficulty in Section \ref{sec:solution}
using the recently developed SOCP relaxation of
OPF.  Fixing any station assignment, the relaxation of OPF is then convex.  Sufficient
conditions are known that guarantee an optimal solution to the nonconvex OPF
problem can be recovered from an optimal solution to its relaxation;
see \cite{Low2014a, Low2014b} for a comprehensive tutorial and references therein.
Even when these conditions are not satisfied, SOCP relaxation is
often exact for practical radial networks, as confirmed also by our simulations.

The second difficulty can be addressed using two different approaches.
The first approach, presented in Section \ref{sec:solution} of this paper,
applies generalized Benders decomposition to the mixed integer convex
relaxation, and is suitable for cases where the distribution network, stations, and EVs
are managed centrally by the same operator.  When
the underlying relaxation of OPF is exact, the generalized Benders decomposition
computes a global optimum.
In Section \ref{sec:numerical results} we illustrate the performance of our centralized solution
through simulations.
The simulation results suggest that the proposed algorithm is effective and computationally efficient for practical application.

In the first approach, the operator needs global information  such as the grid
topology, impedances, operational constraints, background loads, availability
of fully-charged batteries at each station, locations and states of charge of EVs.
The second approach relaxes the binary station assignment variables to real variables in $[0,1]$.
With both relaxations the resulting approximate problem of joint battery swapping scheduling and OPF is a
convex problem.   This allows us to develop distributed solutions that are suitable for cases
where the grid, stations, and EVs are operated by different organizations
that do not share their private information.  Their respective decisions are coordinated
through privacy-preserving information exchanges.
This will be explained in Part II of this paper.

\subsection{Literature}

There is a large literature on EV charging, e.g.,
optimizing charging schedule for various purposes such as demand response,
load profile flattening, or frequency regulation,
e.g., \cite{MaCallawayHiskens2013, gan2013optimal, papadopoulos2013coordination, han2011estimation, o2014rolling};
architecture for mass charging \cite{ChenTong2012}; locational marginal pricing
for EV \cite{li2014distribution}; or the interaction of EV penetration and
the optimal siting and investment of charging stations \cite{Tong2014}.

Sojoudi \emph{et al.} \cite{sojoudi2011optimal} seems to be the first to jointly optimize
EV charging and AC power flow spatially and temporally through semidefinite relaxation.
Zhang \emph{et al.} \cite{zhang2015scalable} extends the joint OPF-charging problem to
multiphase distribution networks and proposes a distributed charging algorithm based on
the alternating direction method of multipliers (ADMM).
Chen \emph{et al.} \cite{chen2014electric} decomposes the joint OPF-charging problem
into an OPF subproblem that is solved centrally by a utility company and a
charging subproblem that is solved in a distributed manner by the EVs coordinated
by a valley-filling signal from the utility.
De Hoog \emph{et al.} \cite{de2015optimal} uses a linear model and formulates EV charging
on a three-phase unbalanced grid as a receding horizon optimization problem.   It shows
that optimizing charging schedule can increase the EV penetration that is sustainable by
the grid from 10--15\% to 80\%.
Linearization is also used in \cite{7031457} to model EV charging on a three-phase
unbalanced grid as a mixed-integer linear program (binary because an EV is either being
charged at peak rate or off).

The literature on battery swapping is comparatively much smaller.
 Tan \emph{et al.} \cite{tan2014queueing} proposes a mixed queueing network
 that consists of a closed queue of batteries and an open queue of EVs to model
 the battery swapping processes, and analyzes its steady-state distribution.
 Yang \emph{et al.} \cite{yang2014dynamic} designs a dynamic operation model
 of a battery swapping station and devises a bidding strategy in power markets.
 %
%
 You \emph{et al.} \cite{7310884} studies the optimal charging schedule of
 a battery swapping station serving electric buses and proposes an efficient
 distributed solution that scales with the number of charging boxes in the station.
 Sarker \emph{et al.} \cite{sarker2015optimal} proposes a day-ahead model for
 the operation of battery swapping stations and uses robust optimization to deal
 with future uncertainty of battery demand and electricity prices.
 Zheng \emph{et al.} in \cite{zheng2014electric} studies the optimal design and
 planning of a battery swapping station in a distribution system to maximize its net present value taking into account life cycle cost of batteries, grid upgrades,
 reliability, operational costs and investment costs.
 Zhang \emph{et al.} \cite{zhang2016benefit} discusses several potential commercial modes of battery swapping and leasing service in China, and presents a benefit analysis from perspectives of utility companies and battery manufacturers.

\section{Problem formulation}\label{sec:system model}

We focus on the scenario where a fleet of EVs and a set of stations\footnote{Throughout this paper stations refer to battery stations.} operate in a region that is supplied by an active distribution network.
We assume the EVs, the stations, and the distribution network are managed centrally by the same operator,
e.g., the State Grid in China.
Periodically, say, every 15 minutes, the system determines a set of EVs that should be scheduled for
battery swapping, e.g., based on their current state of charge or EVs' requests for battery swapping.
At the beginning of the control interval
the system assigns to each EV in the set a station for battery swapping.
The EVs travel to their assigned station to swap their batteries before the end of the current interval,
and batteries returned by the EVs start to be charged at the stations from the next interval.
Our goal is to design an assignment algorithm that
optimizes a weighted sum of electricity generation cost and the distance travelled for battery swapping, while respecting
the operational constraints of the distribution network.

We make two simplifying assumptions.
First we assume that all EVs in the set can arrive at their assigned station and finish battery swapping
before the next interval, so we do not consider scheduling across multiple intervals.
This assumption is reasonable because the geographic area covered by a distribution network is usually
relatively small.  Typically a city substation (50MVA, 110kV) has a service radius of 3--5km, depending on its load density \cite{substationserviceradius}.
Second we ignore the possibility that an EV does not swap its battery as
recommended or swaps its battery at a station different from its assigned station.
These complications affect the initial state at the beginning of the next interval,
but in this paper, we focus only on optimal scheduling in the current interval.

In the following we present a mathematical model of a radial distribution network and formulate our optimal
scheduling problem for battery swapping.
All vectors $x$ in this paper are column vectors; $x^T$ denotes its transpose.

\subsection{Network model}

\begin{figure}
\centering
\includegraphics[width=0.33\textwidth]{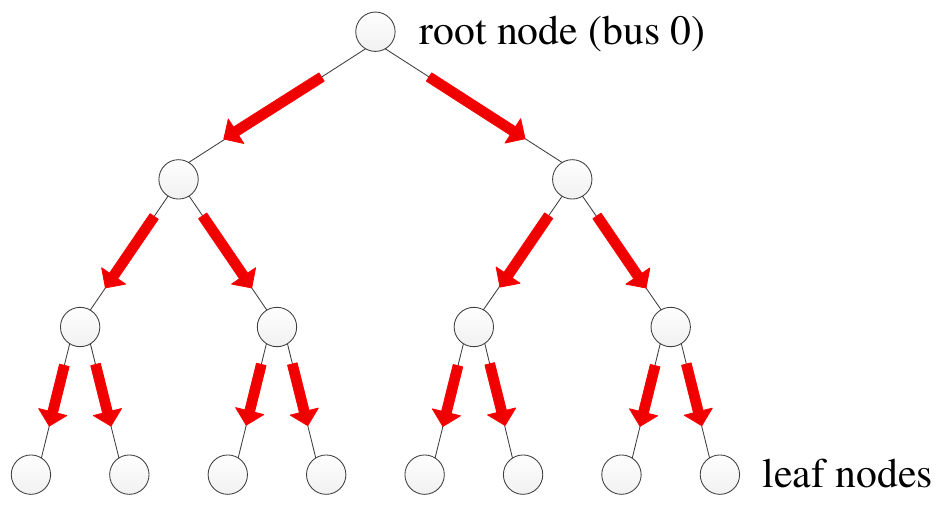}
\caption{Radial topology of $\mathbb{G}$.}
\label{fig:tree}
\end{figure}

Consider a radial distribution network with a connected directed graph
$\mathbb{G}=(\mathbb{N},\mathbb{E})$, where $\mathbb{N}:= \{0,1,2,\dots,N\}$ and $\mathbb{E}\subseteq \mathbb{N} \times \mathbb{N}$. Each node in $\mathbb{N}$ represents a bus and each edge in $\mathbb{E}$ represents a distribution line. We assume $\mathbb{G}$ has a radial (tree) topology with bus 0 representing a substation that extracts power from a transmission network to
feed the loads in $\mathbb G$, as illustrated in Fig.~\ref{fig:tree}.
We orient the graph, without loss of generality, so that each line points away from bus 0.
Denote a line in $\mathbb{E}$ by $(j,k)$ or $j\rightarrow k$ if it points from bus $j$ to bus $k$.
Each bus $j$ (except bus 0) has a unique parent bus $i:=i_j$.
Let $z_{jk}$ be the complex impedance of line $(j,k)\in \mathbb{E}$.  
Let $S_{jk}:=P_{jk}+\textbf{i}Q_{jk}$ denote the \emph{sending-end} complex power
from bus $j$ to bus $k$ where $P_{jk}$ and $Q_{jk}$ denote the real and reactive
power flows.
Let $l_{jk}$ denote the squared magnitude of the complex current from bus $j$ to bus $k$.
Let $v_j$ denote the squared magnitude of the complex voltage phasor of bus $j$.
We assume the voltage $v_0$ of bus 0 is fixed.

Each bus $j$ has a base load $s_j^b := p^b_j + \textbf{i}q^b_j$ (excluding the battery charging loads from stations), where $p_j^b$ and $q_j^b$ denote the real and reactive power.
Each bus $j$ may also have distributed generation $s_j^g:=p_j^g+\textbf{i}q_j^g$.
Let $s_j$ denote the net complex power injection given by
\bqn
s_j & \!\!\!\! :=  \!\!\!\! & \left\{\begin{array}{ll}
		 \!\!\!\! s^g_j - s^b_j - s^e_j  & \text{ if $j$ supplies a station}
		\\
		 \!\!\!\! s^g_j - s^b_j & \text{ otherwise}
		\end{array} \right.
\eqn
where $s^e_j$ denotes the total charging load at bus $j$.
We assume the base loads $s_j^b$ are given and the generations $s_j^g$
and charging loads $s^e_j$ are variables.

We use the \emph{DistFlow} equations proposed by Baran and Wu in \cite{baran1989optimal} 
to model the power flows on the network:
\begin{subequations}\label{eq:df}
\begin{align}\label{eq:df.a}
\sum\limits_{k:(j,k)\in\mathbb{E}}{S_{jk}} &=S_{ij}-z_{ij}l_{ij}+s_j,\quad j\in\mathbb{N}\\\label{eq:df.b}
v_j-v_k &=2 \mathrm{Re} (z_{jk}^HS_{jk})-|z_{jk}|^2l_{jk},\quad j\rightarrow k \in\mathbb{E}\\\label{eq:df.c}
v_jl_{jk} &= |S_{jk}|^2, \quad j\rightarrow k \in\mathbb{E}
\end{align}
\end{subequations}
where $v_j=:|V_j|^2$ and $l_{jk}:=|I_{jk}|^2$.
The equations \eqref{eq:df.a} impose power balance at each bus, \eqref{eq:df.b} model the Ohm's law, and \eqref{eq:df.c} define branch power flows. Note that $S_{i0}:=0$ and $I_{i0}:=0$ if $j=0$ is the substation bus.
When bus $j$ is a leaf node of $\mathbb{G}$, all $S_{jk}=0$ in \eqref{eq:df.a}.
The quantity $z_{ij}|I_{ij}|^2$ is the loss on line $(i,j)$, and hence $S_{ij}-z_{ij}|I_{ij}|^2$ is the \emph{receiving-end} complex power at bus $j$ from $i$.
%
%
%

The complex notation of the \emph{DistFlow} equations \eqref{eq:df} is only a shorthand
for a set of real equations in the real vector variables
$(s, v, l, S) := (p, q, v, l, P, Q) := (p_j, q_j, v_j, l_{jk}, P_{jk}, Q_{jk}, \, j,k \in\mathbb N, \, (j,k)\in\mathbb E)$.
The equations \eqref{eq:df.a}--\eqref{eq:df.b} are linear in these variables
but \eqref{eq:df.c} are quadratic, one of the two sources of nonconvexity in our joint battery swapping
scheduling and OPF problem formulated below.

The operation of the distribution network must meet certain specifications.
The squared voltage magnitudes must satisfy
\begin{subequations}
  \beq\label{eq:con.vol}
  \underline{v}_j \le v_j \le\overline{v}_j, \quad j \in \mathbb{N}
  \eeq
  where $\underline{v}_j$ and $\overline{v}_j$ are given lower and upper bounds on the squared voltage magnitude at bus $j$.
The distributed real and reactive generations must satisfy
  \bq\label{eq:con.gen.p}
  \underline{p}^g_j \le p^g_j \le\overline{p}^g_j, \quad j \in \mathbb{N} \\\label{eq:con.gen.q}
  \underline{q}^g_j \le q^g_j \le\overline{q}^g_j, \quad j \in \mathbb{N}
  \eq
  where $\underline{p}^g_j$, $\overline{p}^g_j$, $\underline{q}^g_j$, and $\overline{q}^g_j$ are given lower and upper bounds on the real and reactive power generation at bus $j$ respectively.
  The power flows on line $(j,k)$ must satisfy
  \beq\label{eq:con.powflow}
  |S_{jk}| \le \overline{S}_{jk}, \quad j\rightarrow k\in\mathbb{E}
  \eeq
 \label{eq:gridconstraints}%
 \end{subequations}
  where $\overline{S}_{jk}$ denotes the capacity of line $(j,k)$.

  The model is quite general.  For example, if a quantity is known and fixed,
  then we set both its upper and lower bounds to the given quantity, e.g., for the
  voltage of the substation bus, $\overline v_0 = \underline v_0$.
  If there is no distributed generation at bus $j$ then
  $\overline p^g_j = \underline p^g_j = \overline q^g_j = \underline q^g_j = 0$.

\subsection{Battery swapping scheduling}

Let $\mathbb{N}_w :=\{1,2,\dots,N_w\} \subseteq \mathbb{N}$ denote the set of buses
that supply electricity to stations.   Their locations are fixed and known.
There is a station connected to each bus $j\in\mathbb N_w$ and we use $j$ to index both
the bus and the station.
The batteries at station $j$ are either charging at a constant rate $r$ or
already fully-charged and ready for swapping.
Denote the total number of batteries and fully-charged batteries at the beginning of the
current control interval by $M_j$ and $m_j$ respectively.

 Let $\mathbb{A}:=\{1,2,\dots,A\}$ denote the set of EVs in the geographic area served by the
distribution network that require battery swapping in the current interval.
Let their states of charge be $(c_a, {a\in\mathbb{A}})$.
Let $u_{aj},~a\in\mathbb{A},~j\in\mathbb{N}_w,$ represent the assignment:
\bqn
u_{aj} & = &\left\{
\begin{split}
& 1,\quad \mathrm{if~EV}~a~\mathrm{is~assigned~to~station}~j\\
& 0,\quad \mathrm{otherwise}
\end{split}
\right.
\eqn
and let ${u}:=(u_{aj}, \, {a\in\mathbb{A}, j\in\mathbb{N}_w}$) denote the vector of assignments.

\vspace{0.055in}
\noindent
\textit{Assumption 1.} $A \leq \sum_{j\in\mathbb{N}_w}{m_j}$.

\vspace{0.055in}
\noindent
Under Assumption 1, there are enough fully-charged batteries in the system for
all EVs in $\mathbb A$ in the current interval.
This can be enforced when choosing the candidate set $\mathbb A$ of EVs
for battery swapping, e.g., ranking EVs according to their states of charge and scheduling in an
increasing order for at most $\sum_{j\in\mathbb{N}_w}{m_j}$ EVs.

The assignment $u$ satisfies the following conditions:
\begin{subequations}
  \bq\label{eq:con.tax1}
  \sum\limits_{j\in\mathbb{N}_w}{u_{aj}} & = & 1,\quad a\in\mathbb{A} \\\label{eq:con.tax2}
  \sum\limits_{a\in\mathbb{A}}{u_{aj}} & \le & m_j,\quad j\in\mathbb{N}_w
  \eq
i.e., exactly one station is assigned to every EV and every assigned station has enough
fully-charged batteries to serve EVs.

The system knows the current location of every EV $a$ and therefore can calculate
the distance $d_{aj}$ from its current location to the assigned station $j$.  If the EV
is not currently carrying passengers and can go to swap its battery immediately, then $d_{aj}$ is the travel distance from its current location to station $j$, e.g.,
calculated from a routing application (such as Google map).  If the EV must
first complete its current passenger run before going to station $j$, then the distance
$d_{aj}$ is the travel distance from its current location to the destination of
its passengers and then to station $j$.   The assigned station $j$ must be within each EV $a$'s driving range, i.e.,
      \beq\label{eq:con.tax3}
      u_{aj}d_{aj} \le  \gamma_a c_a,\quad j\in\mathbb{N}_w, a\in\mathbb{A}
      \eeq
\label{eq:uconstraints}%
\end{subequations}
where $c_a$ is EV $a$'s current state of charge and $\gamma_a$ is its driving range per unit state of charge.


Since every EV produces a depleted battery that needs to be charged
at rate $r$, we can express the net power injection $s_j=p_j+\textbf{i}q_j$ at
bus $j$ in terms of assignment $u$ as:
\begin{subequations}
\bq
\label{eq:inj.p}
p_j & \!\!\!\!\! = \!\!\!\!\! & \left\{
\begin{split}
&p_j^g-p_j^b-r\!\left(M_j-m_j+\sum_{a\in\mathbb{A}}{u_{aj}}\right)\!,& \!\!\! & j\in\mathbb{N}_w\\
&p_j^g-p_j^b, & \!\!\! & j\in\mathbb{N}\setminus\mathbb{N}_w
\end{split}
\right.
\\
\label{eq:inj.q}
q_j & \!\!\!\! = \!\!\! & \ \ q_j^g-q_j^b,\qquad j\in\mathbb{N}
\eq
\label{eq:inj}%
\end{subequations}

Let $f_j:\mathbb{R}\rightarrow \mathbb{R}$ models the generation cost at bus $j$, e.g., for
a distributed gas generator.   We assume all $f_j$ are strictly convex increasing functions
\cite{sojoudi2011optimal, chen2014electric, zhang2015scalable}.
We are interested in the following optimization problem:
\bq\label{eq:primalprob}
\min\limits_{\bm{u},\bm{s}, \bm{s}^g, \atop \bm{v},\bm{l}, \bm{S}} &&
	\sum\limits_{j\in\mathbb{N}}{f_j(p_j^g)} + \alpha
	\sum\limits_{a\in\mathbb{A}}\sum\limits_{j\in\mathbb{N}_w} d_{aj} u_{aj}
\\\nonumber
\rm{s.t.} && \eqref{eq:df} \eqref{eq:gridconstraints} \eqref{eq:uconstraints} \eqref{eq:inj},
	\ u_{aj}\in \{0, 1\}
\eq
where $\sum\limits_{a\in\mathbb{A}}\sum\limits_{j\in\mathbb{N}_w} u_{aj}d_{aj}$ is the
total travel distance of EVs and $\alpha>0$ is a weight that makes the generation cost and the travel
distance comparable.

\section{Solution}\label{sec:solution}

The joint battery swapping scheduling and OPF problem \eqref{eq:primalprob} is generally difficult to solve
because \eqref{eq:df.c} is nonconvex, as mentioned above, and $u$ is discrete.
Our solution strategy has two steps.

\vspace{0.1in}
\noindent
\textbf{1. SOCP relaxation.}
We first relax the nonconvex constraint \eqref{eq:df.c} into a second-order cone.
Specifically, replace the DistFlow equations \eqref{eq:df} by
\begin{subequations}\label{eq:df2}
\begin{align}\label{eq:df2.a}
\sum\limits_{k:(j,k)\in\mathbb{E}}{S_{jk}} &=S_{ij}-z_{ij}l_{ij}+s_j,\quad j\in\mathbb{N}\\\label{eq:df2.b}
v_j-v_k &=2 \mathrm{Re} (z_{jk}^HS_{jk})-|z_{jk}|^2l_{jk},\quad j\rightarrow k \in\mathbb{E}\\\label{eq:df2.c}
v_jl_{jk} &\geq |S_{jk}|^2, \quad j\rightarrow k \in\mathbb{E}
\end{align}
\end{subequations}
Then the SOCP relaxation of the problem \eqref{eq:primalprob}
is:
\bq\label{eq:primalprob2}
\min\limits_{\bm{u}, \bm{s}, \bm{s}^g, \atop \bm{v},\bm{l}, \bm{S}} &&
	\sum\limits_{j\in\mathbb{N}}{f_j(p_j^g)} + \alpha
	\sum\limits_{a\in\mathbb{A}}\sum\limits_{j\in\mathbb{N}_w} d_{aj} u_{aj}
\\\nonumber
\rm{s.t.} && \eqref{eq:df2} \eqref{eq:gridconstraints} \eqref{eq:uconstraints} \eqref{eq:inj},
	\ u_{aj}\in \{0, 1\}
\eq
Fix any assignment $u\in\{0,1\}^A$.  Then the problem \eqref{eq:primalprob2} is a convex
problem.   It is a relaxation of the problem \eqref{eq:primalprob}, given $u$, in the sense that the
optimal objective value of the relaxation \eqref{eq:primalprob2} lower bounds that of the original problem
\eqref{eq:primalprob}.   If an optimal solution to the relaxation \eqref{eq:primalprob2}
attains equality in \eqref{eq:df2.c} then the solution is also feasible, and therefore
\emph{optimal}, for the original problem \eqref{eq:primalprob}.  In this case, we say
that the SOCP relaxation is \emph{exact}.
Sufficient conditions are known that guarantee the exactness of the SOCP relaxation;
see \cite{Low2014a, Low2014b} for a comprehensive tutorial and references therein.
Even when these conditions are not satisfied, SOCP relaxation for practical radial networks
is often exact, as confirmed also by our simulations in
Section \ref{sec:numerical results}.

Hence we will solve \eqref{eq:primalprob2} instead of \eqref{eq:primalprob}.

\vspace{0.1in}
\noindent
\textbf{2. Generalized Benders decomposition.}
To deal with the discrete variables in \eqref{eq:primalprob2}, we
apply generalized Benders decomposition.
Benders decomposition was first proposed in \cite{benders1962partitioning}
for problems where, when a subset of the variables are fixed, the remaining
subproblem is a linear program.  It is extended in \cite{geoffrion1972generalized}
to the broader class of problems where the remaining subproblem is a convex
program.   We now apply it to solving \eqref{eq:primalprob2}.

Denote the continuous variables by $x := (s, {s}^g,  {v},{l}, {S})$
and the discrete variables by $u$.   Denote the objective function by
\bqn
F(x, u) & := & 	\sum\limits_{j\in\mathbb{N}}{f_j(p_j^g)} + \alpha
	\sum\limits_{a\in\mathbb{A}}\sum\limits_{j\in\mathbb{N}_w} d_{aj} u_{aj}
\eqn
Given any $u$, $F(x, u)$ is convex in $x$ since $f_j$'s are assumed to be
strictly convex.   Denote the constraint set for $x$ by
\bqn
\mathbb X & := & \{ x\in \mathbb R^{(5|\mathbb N|+3|\mathbb E|)} \, :\,
	x \text{ satisfies } \eqref{eq:gridconstraints} \eqref{eq:df2}  \}
\eqn
the constraint set for $u$ by
\bqn
\mathbb U & := & \{ u\in \{0, 1\}^{A N_w} \, :\,
	u \text{ satisfies } \eqref{eq:uconstraints} \}
\eqn
and the constraints \eqref{eq:inj} on $(x, u)$ by $G(x, u) = 0$.
Then the relaxation \eqref{eq:primalprob2} takes the standard form
for generalized Benders decomposition:
\bq\label{eq:MINLP}
\min\limits_{x,u} && F(x,u)\\\nonumber
\rm{s.t.} && G(x,u) = 0,~x\in\mathbb{X},~u\in\mathbb{U}
\eq
where
$F:\mathbb{R}^{(5|\mathbb N|+3|\mathbb E|)}\times \{0, 1\}^{AN_w}\rightarrow \mathbb{R}$
is a scalar-valued function, and
$G:\mathbb{R}^{(5|\mathbb N|+3|\mathbb E|)}\times \{0, 1\}^{AN_w}\rightarrow \mathbb{R}^{2|\mathbb N|}$
is a vector-valued constraint function.
Fixing any $u\in \mathbb U$, \eqref{eq:MINLP} is a convex subproblem in $x$.
We now apply generalized Benders decomposition of \cite{geoffrion1972generalized}
to \eqref{eq:MINLP}.

Write \eqref{eq:MINLP} in the following equivalent form:
\begin{subequations}
\bq \label{eq:projection}
\min\limits_{u} \ \ W(u) & &
\mathrm{s.t.} \qquad u\in \mathbb{U} \cap \mathbb{W}
\eq
where
\beq\label{eq:Wh}
\begin{split}
W(u) := \ &\min\limits_{x\in\mathbb X} \quad F(x,u) \\
 & \ \mathrm{s.t.} \quad \ G(x,u)= 0
\end{split}
\eeq
and
\beq\label{eq:W}
\mathbb{W}:=\{u:G(x,u) = 0 \mathrm{~for~some}~x \in\mathbb{X}\}
\eeq
\label{eq:proj}
\end{subequations}
The problem \eqref{eq:Wh}, called the slave problem,
 is convex and much easier to solve than \eqref{eq:MINLP}.
The set $\mathbb{W}$ consists of all $u$ for which \eqref{eq:Wh} is feasible
and hence $\mathbb{U} \cap \mathbb{W}$ is the projection of the feasible region
of \eqref{eq:MINLP} onto the $u$-space.
The central idea of generalized Benders decomposition is to invoke the dual
representations of $W(u)$ and $\mathbb{W}$ to derive the following equivalent problem
to \eqref{eq:proj} (see \cite[Theorems 2.2 and 2.3]{geoffrion1972generalized}):
\bqn
\min\limits_{u\in\mathbb{U}} && \sup\limits_{\mu\in \mathbb{R}^{2|\mathbb N|}}
	\left\{\min\limits_{x\in\mathbb{X}}\left\{F(x,u)+\mu^T G(x,u)\right\}\right\}\\\nonumber
\mathrm{s.t.} && \min\limits_{x\in\mathbb{X}} \left\{\lambda^T G(x,u)\right\} = 0,
	\quad \forall \lambda \in \mathbb{R}^{2|\mathbb N|}
\eqn
Here $\lambda$ and $\mu$ are Lagrangian multiplier vectors
for $\mathbb W$ and $W(u)$ respectively.
This problem is equivalent to:
\bq\label{eq:primalmaster}
\min \limits_{u\in\mathbb{U},u_0}&& u_0\\\nonumber
\mathrm{s.t.} && u_0 \ge \min\limits_{x\in\mathbb{X}}\left\{F(x,u)+\mu^T G(x,u)\right\},
		\ \ \forall \mu \in \mathbb{R}^{2|\mathbb N|}
		\\\nonumber
&& \min\limits_{x\in\mathbb{X}} \left\{\lambda^T G(x,u)\right\} = 0,
		\quad \forall \lambda \in \mathbb{R}^{2|\mathbb N|}
\eq
In summary, the series of manipulations has transformed the relaxation
\eqref{eq:primalprob2} into the master problem \eqref{eq:primalmaster}.

Since \eqref{eq:primalmaster} has uncountably many constraints with all possible $\lambda$'s and $\mu$'s, it is neither practical nor necessary to enumerate all constraints in solving \eqref{eq:primalmaster}.
Generalized Benders decomposition starts by solving a relaxed version of
\eqref{eq:primalmaster} that ignores all but a few constraints.  If a solution of
the relaxed version of \eqref{eq:primalmaster} satisfies all the ignored constraints, then
it is an optimal solution of \eqref{eq:primalmaster} and the algorithm
terminates.   Otherwise, the solution process of the relaxed version of \eqref{eq:primalmaster}
will identify one $\mu$ or $\lambda$ for which the constraints are violated.
These constraints are then added to the relaxed version of \eqref{eq:primalmaster}, and
the cycle repeats.

Specifically the Benders decomposition algorithm for \eqref{eq:primalprob2}
(or equivalently \eqref{eq:MINLP}) is as follows.
\begin{itemize}
  \item \textbf{Step 1.}  Pick any $\bar u\in \mathbb{U} \cap \mathbb{W}$.
  Solve \eqref{eq:Wh} with  $u = \bar u$ to obtain an optimal Lagrangian multiplier
  vector $\bar \mu$. Let $n_{\mu}=1$, $n_{\lambda}=0$, $\mu^1=\bar \mu$,
  and $UBD=W(\bar u)$, where $n_{\mu}$, $n_{\lambda}$ are counters for the two
  types of constraints in \eqref{eq:primalmaster}, and
  $UBD$ denotes an upper bound on the optimal value of \eqref{eq:MINLP}.

  \item \textbf{Step 2.} Solve the current relaxed master problem:
  \bq\label{eq:relaxedmaster}
\min \limits_{u\in\mathbb{U},u_0}&& u_0\\\nonumber
\mathrm{s.t.} && u_0 \ge \min\limits_{x\in\mathbb{X}}\left\{F(x,u)+\left({\mu^i}\right)^T G(x,u)\right\}, \\\nonumber
&& \qquad\qquad\qquad\qquad\qquad\ \   i=1,\dots,n_{\mu} \\\nonumber
&& \min\limits_{x\in\mathbb{X}} \left\{\left({\lambda^i}\right)^T G(x,u)\right\} = 0,\\\nonumber
&& \qquad\qquad\qquad\qquad\qquad\ \   i=1,\dots,n_{\lambda}
\eq

Let $(\hat u,\hat u_0)$ be the optimal solution to \eqref{eq:relaxedmaster}.
Clearly $\hat u_0$ is a lower bound on the optimal value of \eqref{eq:MINLP} since
the constraints in \eqref{eq:primalmaster} are relaxed to a smaller set of constraints in
\eqref{eq:relaxedmaster}.
Terminate the algorithm if $UBD-\hat u_0 \le \epsilon$, where $\epsilon>0$ is
a sufficiently small threshold.

  \item \textbf{Step 3.} Solve the dual problem of \eqref{eq:Wh} with $u = \hat u$.
  The solution falls into the following two cases.

  \bee
  \item \textbf{Step 3a.} \emph{The dual problem of \eqref{eq:Wh} has a finite solution $\hat \mu$.}
  $W(\hat u)$ is finite. Let $UBD=\min\{UBD,W(\hat u)\}$.
  Terminate the algorithm if $UBD-\hat u_0 \le \epsilon$.
  Otherwise, increase $n_{\mu}$ by 1 and let $\mu^{n_{\mu}}=\hat \mu$. Return to \textbf{Step 2}.
  \item \textbf{Step 3b.} \emph{The dual problem of \eqref{eq:Wh} has an unbounded solution.}
  Then \eqref{eq:Wh} is infeasible. Determine $\hat \lambda$ through a feasibility check
  problem and its dual \cite{7287792}. Increase $n_{\lambda}$ by 1 and let
  $\lambda^{n_{\lambda}}=\hat \lambda$. Return to \textbf{Step 2}.
  \eee

\end{itemize}

We make three remarks.  First, the slave problem \eqref{eq:Wh} is convex and hence can generally
be solved efficiently.  The relaxed problem \eqref{eq:relaxedmaster} involves discrete variables
and are generally nonconvex, but it is much simpler than the original problem \eqref{eq:MINLP}.
Second, for our problem, \eqref{eq:relaxedmaster} turns out to be a binary linear program because
both $F$ and $G$ are separable functions in $(x,u)$ of the form:
\bqn
F(x, u) & =: & f_1(x) + f_2(u)
\\
G(x, u) & =: & g_1(x) + g_2(u)
\eqn
where $f_2$ and $g_2$ are both linear in $u$.
Indeed the constraints in \eqref{eq:relaxedmaster} are
\bqn
u_0 - f_2(u) - \left({\mu^i}\right)^Tg_2(u)  & \ge &
\min\limits_{x\in\mathbb{X}}\left\{f_1(x)+\left({\mu^i}\right)^T g_1(x)\right\}
\\
&& \qquad\qquad\qquad   i=1,\dots,n_{\mu}
\\
\left({\lambda^i}\right)^T g_2(u) & = & - \min\limits_{x\in\mathbb{X}} \left({\lambda^i}\right)^T g_1(x)
\\
&& \qquad\qquad\qquad   i=1,\dots,n_{\lambda}
\eqn
where the left-hand side is linear in the variable $u$ and the right-hand side is independent of $u$.
Hence, in each iteration, the algorithm solves a binary linear program \eqref{eq:relaxedmaster}
and a convex program \eqref{eq:Wh}.
Finally, every time Step 2 is entered, one or two additional constraints are added to the
binary linear program \eqref{eq:relaxedmaster}.  This generally makes \eqref{eq:relaxedmaster}
harder to compute but also a better approximation of \eqref{eq:primalmaster}.
It is proved in \cite[Theorem 2.4]{geoffrion1972generalized}
that the algorithm will terminate in finite steps since $\mathbb U$ is discrete
and finite.

\section{Numerical results}\label{sec:numerical results}

In this section, we evaluate the proposed algorithm through numerical simulations
using a 56-bus distribution feeder of Southern California Edison (SCE) with a radial structure.
More details about the feeder can be found in \cite[Figure 2, TABLE I]{farivar2012optimal}.
We add 4 distributed generators and 4 stations at different buses. The setup of the distributed generators is given in Table \ref{tab:dgsetup}\footnote{The units of the real power, reactive power, cost (for the whole control interval), distance and weight in this paper are MW, Mvar, \$, km and \$/km, respectively.}.
The 4 stations are assumed to be uniformly located in a 4km$\times$4km square area supplied by the distribution feeder, as shown in Table \ref{tab:stationsetup}.
Suppose in a certain control interval, there are $A$ number of EVs that request battery swapping.
Their current locations are randomized uniformly within the square area
while their destinations are ignored. We use the Euclidean distance for $d_{aj}$.
For convenience, it is assumed that $M_j=m_j,~m_j=A,~j\in\mathbb{N}_w$, which means in each station batteries are all fully-charged and sufficient to serve all EVs.
We assume all EVs have sufficient battery energy to reach any of the 4 stations, which means \eqref{eq:con.tax3} is readily satisfied. The extension to the general case where each EV has a limited driving range and can only reach some of the stations is straightforward.
The constant charging rate is $r=0.01$MW \cite{yilmaz2013review} at all stations. To make the two components of the objective comparable, we set the weight $\alpha$ to be 0.02\$/km first \cite{eia2011annual}, and will then allow it to take different values to reveal its impact. Note that due to the randomness of EVs' initialized locations, we conduct 10 simulation runs for each case setup. All numerical tests are run on a laptop with Intel Core i7-3632QM CPU@2.20GHz, 8GB RAM, and 64-bit Windows 10 OS.

%

\begin{table}[htbp]
  \tiny
  \caption{Setup}\label{tab:setup}
  \begin{center}
  \subtable[Distributed generator]{
  \begin{tabular}{ccccccc}
  \toprule
     Bus & $\overline{p}_j^g $ & $\underline{p}_j^g $ & $\overline{q}_j^g $  & $\underline{q}_j^g $  & Cost function \\
  \midrule
  \rowcolor[gray]{.9}   1 & 4 & 0 & 2 & -2 & $0.3 {p^g}^2 + 30 p^g$ \\
    4 & 2.5 & 0 & 1.5 & -1.5 & $0.1 {p^g}^2 + 20 p^g$ \\
  \rowcolor[gray]{.9}   26 & 2.5 & 0 & 1.5 & -1.5 & $0.1 {p^g}^2 + 20 p^g$ \\
    34 & 2.5 & 0 & 1.5 & -1.5 & $0.1 {p^g}^2 + 20 p^g$ \\
  \bottomrule
  \end{tabular}\label{tab:dgsetup}
  }
  \subtable[Station]{
  \begin{tabular}{ccccccc}
  \toprule
     Bus &  Location & $M_j$ & $m_j$   \\
  \midrule
  \rowcolor[gray]{.9}   5 & (1,1) & $m_j$ & $ A$  \\
    16 & (3,1) & $m_j$ & $ A$ \\
  \rowcolor[gray]{.9}  31 & (1,3) & $m_j$ & $A$  \\
    43 & (3,3) & $m_j$ & $A$  \\
  \bottomrule
  \end{tabular}\label{tab:stationsetup}
  }
  \end{center}
\end{table}

\begin{figure}
\centering
\subfigure[]{
	\includegraphics[width=0.21\textwidth]{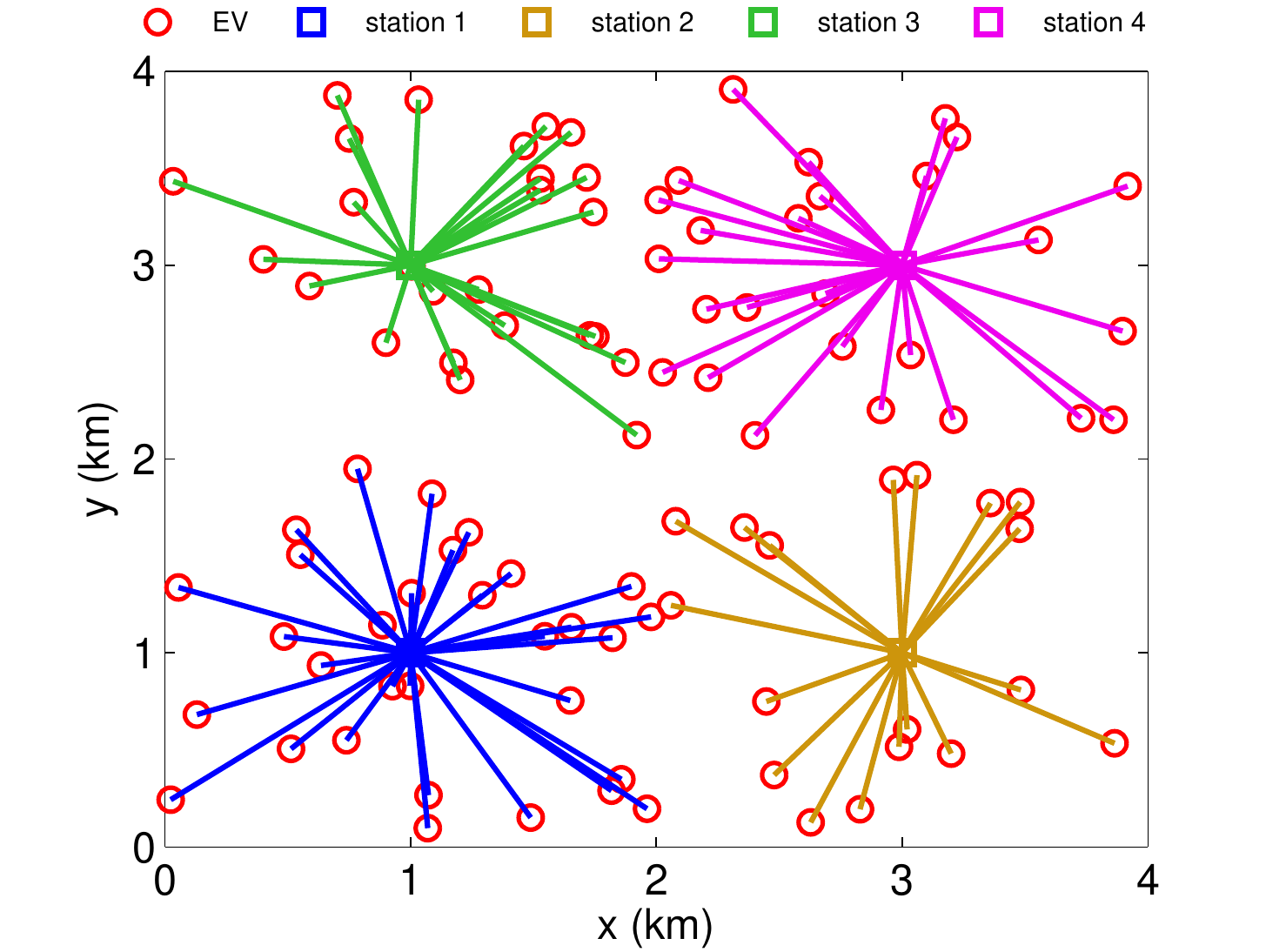}
	\label{fig:nrt400}
}
\hspace{-0.3in}
\subfigure[]{
	\includegraphics[width=0.21\textwidth]{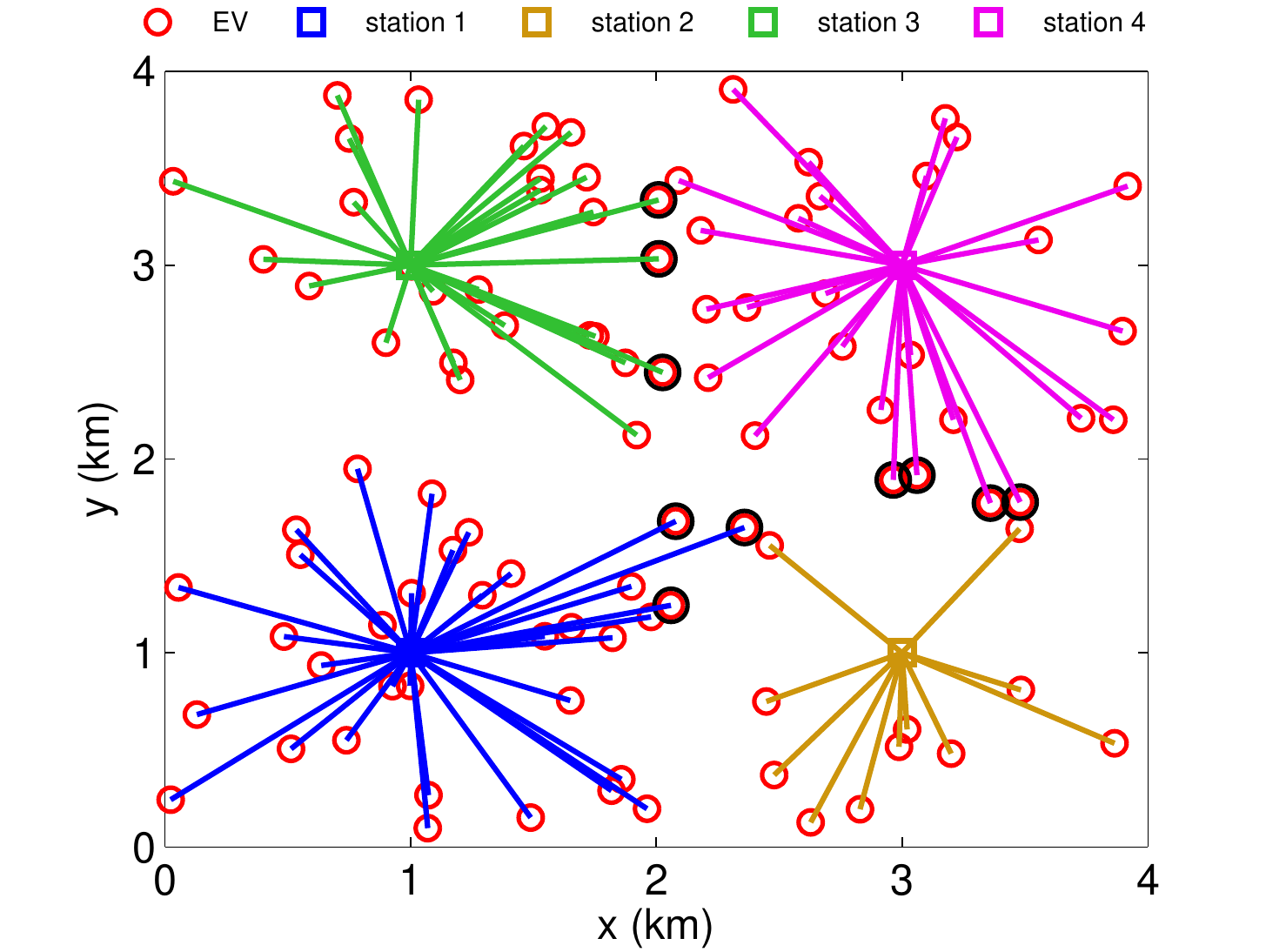}
	\label{fig:opt400}
}
\caption{\#EVs=100 (a) Nearest-station policy. (b) Optimal assignment.}
\end{figure}

\begin{figure}
\centering
\subfigure[]{
	\includegraphics[width=0.21\textwidth]{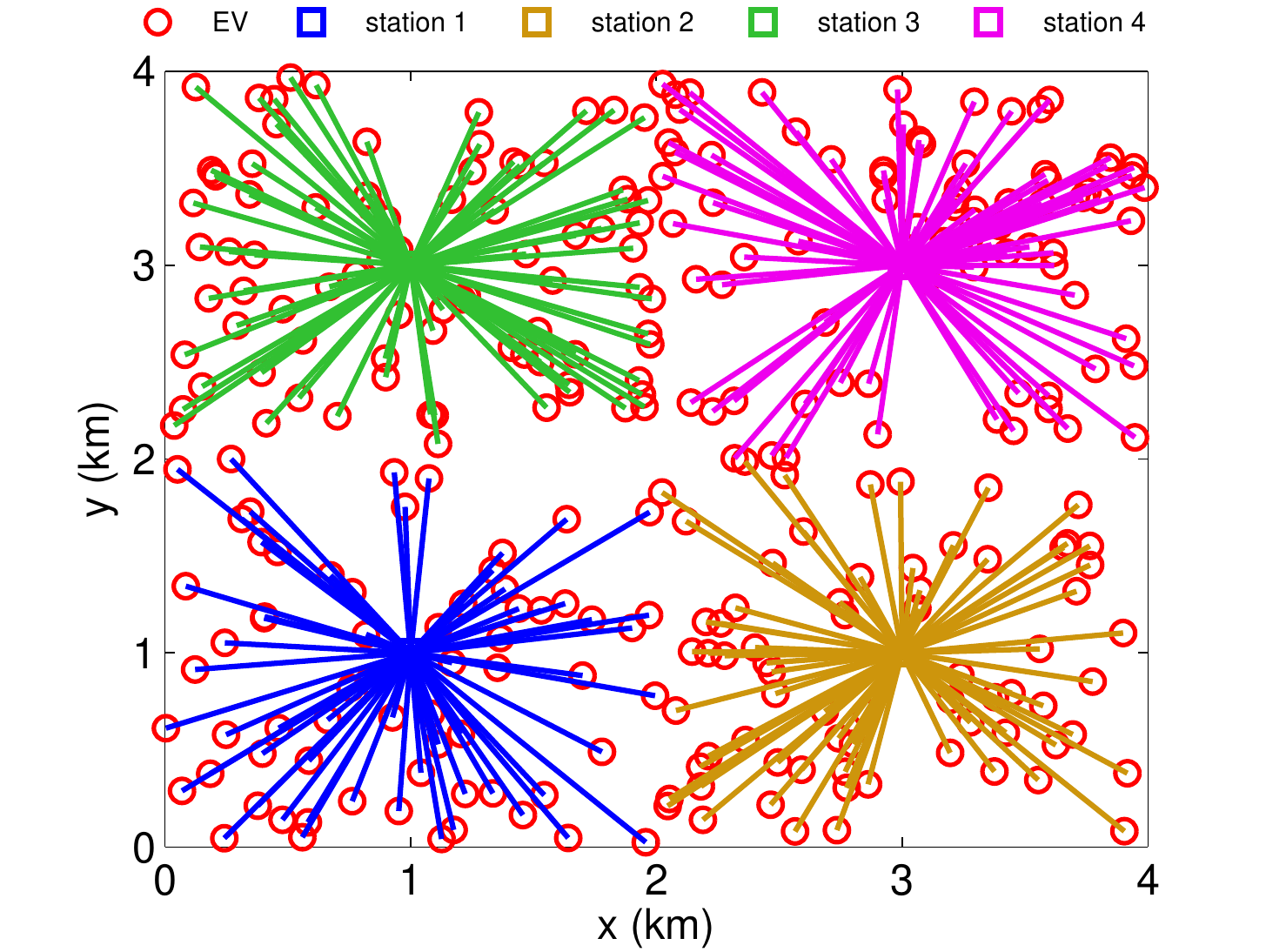}
	\label{fig:nrt800}
}
\hspace{-0.3in}
\subfigure[]{
	\includegraphics[width=0.21\textwidth]{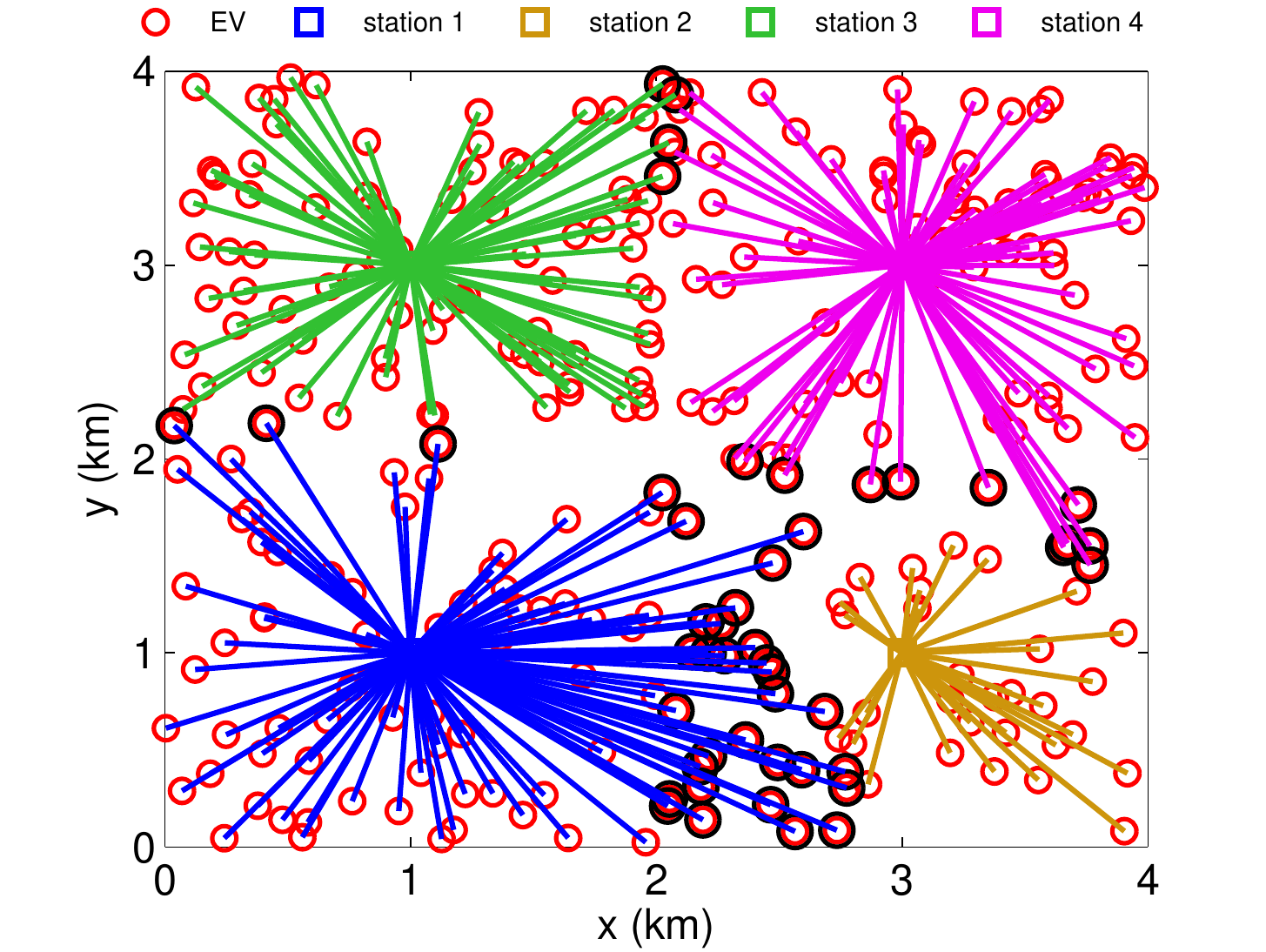}
	\label{fig:opt800}
}
\caption{\#EVs=300 (a) Nearest-station policy. (b) Optimal assignment.}
\end{figure}

\noindent
\textbf{Nearest-station policy.}
Without optimization, the default policy is that all EVs head for their nearest stations to swap batteries.
This is shown in Fig.~\ref{fig:nrt400} and Fig.~\ref{fig:nrt800} for two specific cases with 100 and 300 EVs, respectively.
In practice this myopic policy can lead to battery shortage at a station if many EVs cluster around that station
due to correlation in traffic patterns.
Moreover it can cause voltage instability: the voltage magnitudes of some buses drop below the threshold 0.95 p.u. in the 300-EV case, as shown in Table \ref{tab:opf.bus.withoutsche}.

\begin{table}[htbp]
  \tiny
  \caption{Partial Bus data under nearest-station policy (300 EVs)}\label{tab:opf.bus.withoutsche}
  \begin{center}
  \begin{tabular}{ccccc}
  \toprule
    Bus & $|V_j|$ (p.u.) & $p_j^g $ & $q_j^g $  & $r\sum\limits_{a\in\mathbb{A}}{u_{aj}}$   \\
  \midrule
  \rowcolor[gray]{.9}  1 & 1.050 & 0.571 & 0.000 & / \\
   4 & 1.047 & 2.500 & 0.663 & / \\
  \rowcolor[gray]{.9} 5  & 1.031 & / & / & 0.660 \\
   16 & 0.941 & / & / & 0.700 \\
  \rowcolor[gray]{.9}     18 & 0.948 & / & / & / \\
  19 & 0.944 &  / & / & /\\
  \rowcolor[gray]{.9}     26 & 1.050 & 2.500 & 0.410 & / \\
  31  & 1.020 & / & / & 0.830 \\
  \rowcolor[gray]{.9}     34 & 1.044 & 2.500 & 1.500 & / \\
   43  & 1.015 & / & / & 0.810 \\
  \bottomrule
  \end{tabular}
  \end{center}
\end{table}

\noindent
\textbf{Optimal assignment.}
Fig.~\ref{fig:opt400} and Fig.~\ref{fig:opt800} show the optimal assignments computed using the proposed algorithm for the above two cases, respectively.
The nearest stations are not assigned to some of the EVs (marked black in the figures) when grid operational constraints such as voltage stability are taken into account.
The numbers of such EVs is higher in the 300-EV case than that in the 100-EV case.
The tradeoff between EVs' total
travel distance and the total generation cost is optimized.
For comparison with Table \ref{tab:opf.bus.withoutsche}, the corresponding partial OPF results of the 300-EV case are listed in Table \ref{tab:opf.bus}.
As we can see from Table \ref{tab:opf.bus}, the outputs (2.500 MW) of the distributed generators
at buses 4, 26 and 34 have reached their full capacity (2.5 MW) while the injection (0.520 MW) at bus 1 (root bus) is far from its capacity (4 MW).  This is consistent with our intuition that distributed generators that are closer to users
and potentially cheaper than power from the transmission grid are favored in OPF.
Under the optimal assignment, the deviations of voltages from their nominal value are all less
than the 5\%.


\begin{table}[htbp]
  \tiny
  \caption{Partial Bus data under optimal assignment (300 EVs)}\label{tab:opf.bus}
  \begin{center}
  \begin{tabular}{ccccc}
  \toprule
    Bus & $|V_j|$ (p.u.) & $p_j^g $ & $q_j^g $  & $r\sum\limits_{a\in\mathbb{A}}{u_{aj}}$   \\
  \midrule
  \rowcolor[gray]{.9}  1 & 1.050 & 0.520 & 0.000 & / \\
   4 & 1.048 & 2.500 & 0.590 & / \\
  \rowcolor[gray]{.9} 5  & 1.025 & / & / & 0.990 \\
  15  & 0.981 & / & / & / \\
  \rowcolor[gray]{.9}   16 & 0.974 & / & / & 0.300 \\
  17  & 0.980 & / & / & / \\
  \rowcolor[gray]{.9}     18 & 0.973 & / & / & / \\
  19 & 0.969 &  / & / & /\\
  \rowcolor[gray]{.9}     26 & 1.050 & 2.500 & 0.439 & / \\
  31  & 1.019 & / & / & 0.840 \\
  \rowcolor[gray]{.9}     34 & 1.044 & 2.500 & 1.500 & / \\
  43  & 1.013 & / & / & 0.870 \\
  \bottomrule
  \end{tabular}
  \end{center}
\end{table}

\begin{figure}
\centering
\subfigure[]{
	\includegraphics[width=0.21\textwidth]{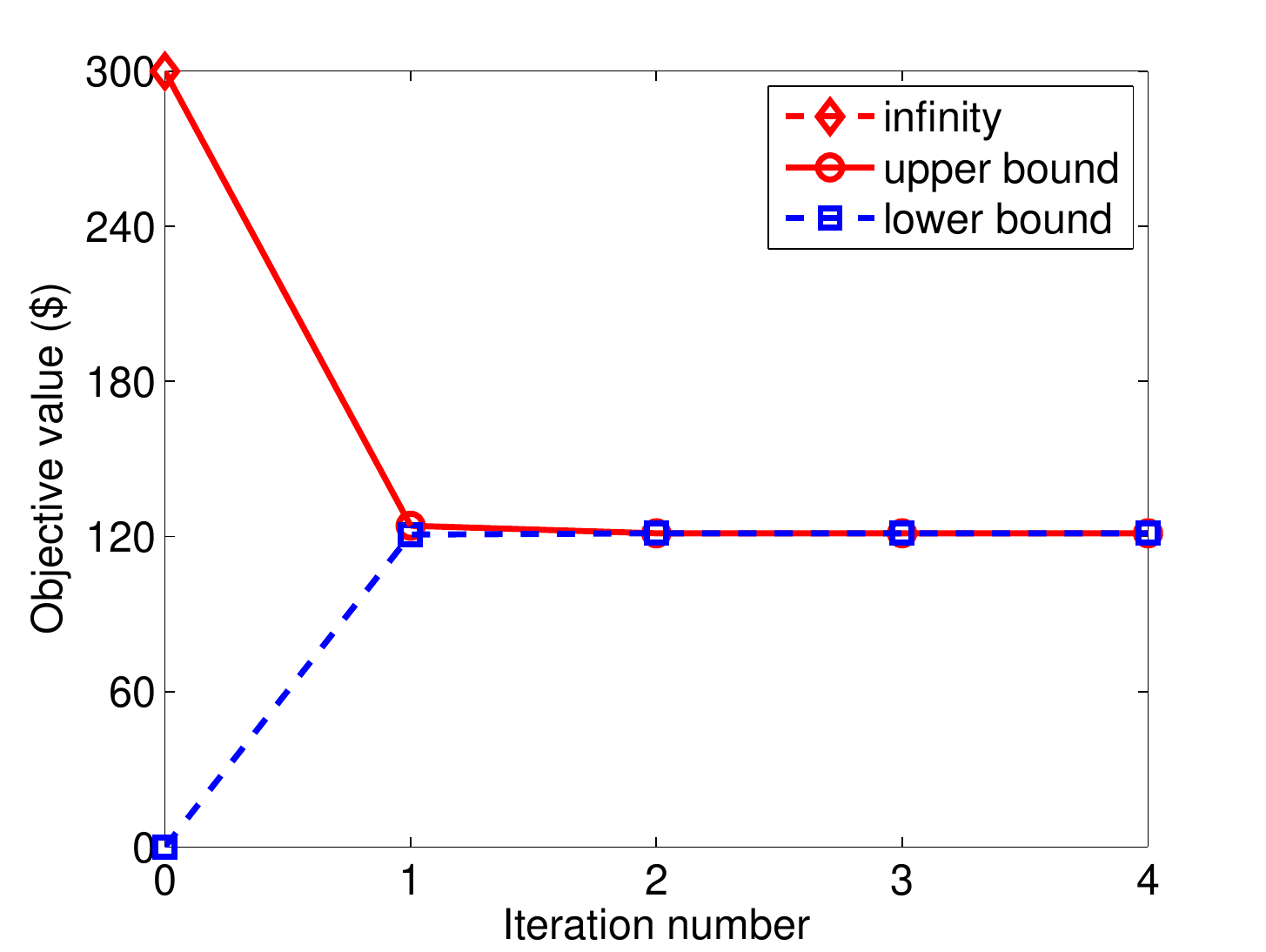}
	\label{fig:conv400}
}
\hspace{-0.3in}
\subfigure[]{
	\includegraphics[width=0.21\textwidth]{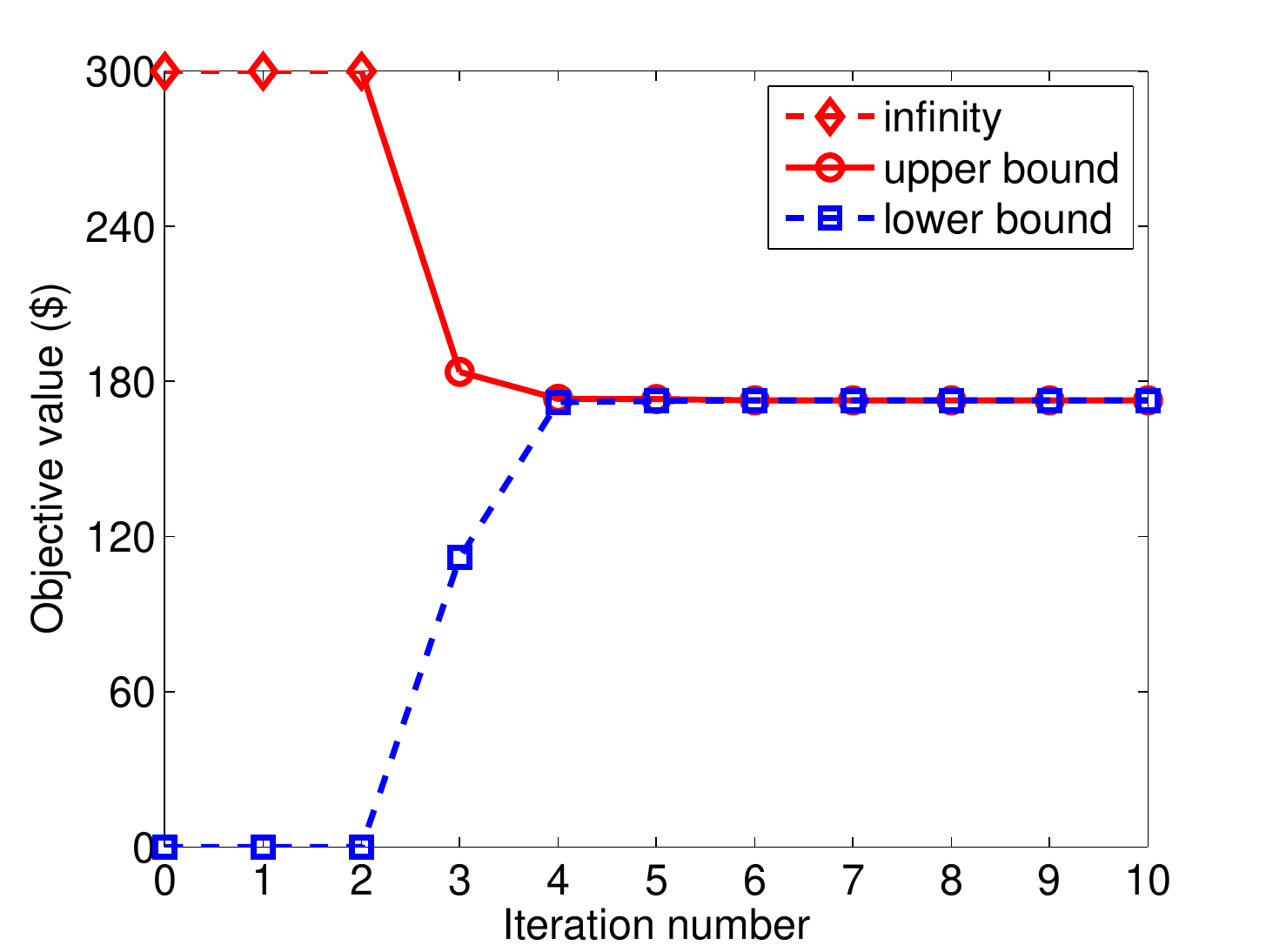}
	\label{fig:conv800}
}
\caption{Convergence of generalized Benders decomposition (a) \#EVs=100. (b) \#EVs=300.}
\end{figure}

\noindent
\textbf{Optimality of generalized Benders decomposition.}
The upper and lower bounds on the optimal objective values for the above two cases are plotted in Fig.~\ref{fig:conv400} and Fig.~\ref{fig:conv800}, respectively, as the algorithm iterates between the master and slave problems. Basically, more iterations are required for larger-scale cases since it usually takes more iterations to attain an initial feasible solution. However, once we have a feasible solution, the gap between the upper and lower bounds starts to shrink rapidly and the convergence to optimality is achieved within a few iterations.

\noindent
\textbf{Exactness of SOCP relaxation.}
We check whether the solution computed by generalized Benders decomposition attains equality
in \eqref{eq:df2.c}, i.e., whether the solution satisfies power flow equations and is implementable.
Our final result confirms the exactness of the SOCP relaxation for the above two cases,
and the relaxation is exact for most other cases we have tested on. Due to space limit, only
some partial data of the 300-EV case are shown in Table \ref{tab:exactness}.

In summary, SOCP relaxation and generalized Benders decomposition have solved our joint battery swapping scheduling and OPF problem \eqref{eq:primalprob} exactly.

\begin{table}[htbp]
  \tiny
  \caption{Exactness of SOCP relaxation (partial results for 300 EVs)}\label{tab:exactness}
  \begin{center}
  \begin{tabular}{ccccc}
  \toprule
  \multicolumn{2}{c}{Bus} & \multirow{2}{*}{ $v_jl_{jk} $} & \multirow{2}{*}{ $|S_{jk}|^2$} & \multirow{2}{*}{Residual}  \\
   From & To & & & \\
  \midrule
  \rowcolor[gray]{.9}   1 & 2 &  0.271 &  0.271 & 0.000 \\
    2 & 3 &  0.006 &  0.006 & 0.000 \\
  \rowcolor[gray]{.9}  2 & 4 & 0.202 & 0.202 & 0.000 \\
    4 & 5 & 1.369 & 1.369 & 0.000 \\
  \rowcolor[gray]{.9}  4 & 6 &  0.005 &  0.005 & 0.000 \\
      4 & 7 & 1.952 & 1.952 & 0.000 \\
  \rowcolor[gray]{.9}  7 & 8 &  1.691 &  1.691 & 0.000 \\
      8 & 9 &  0.009 &  0.009  & 0.000 \\
  \rowcolor[gray]{.9}  8 & 10 & 1.269 & 1.269 & 0.000 \\
        10 & 11 & 1.092 & 1.092 & 0.000 \\
  \bottomrule
  \end{tabular}
  \end{center}
\end{table}

\begin{figure}
\centering
	\includegraphics[width=0.23\textwidth]{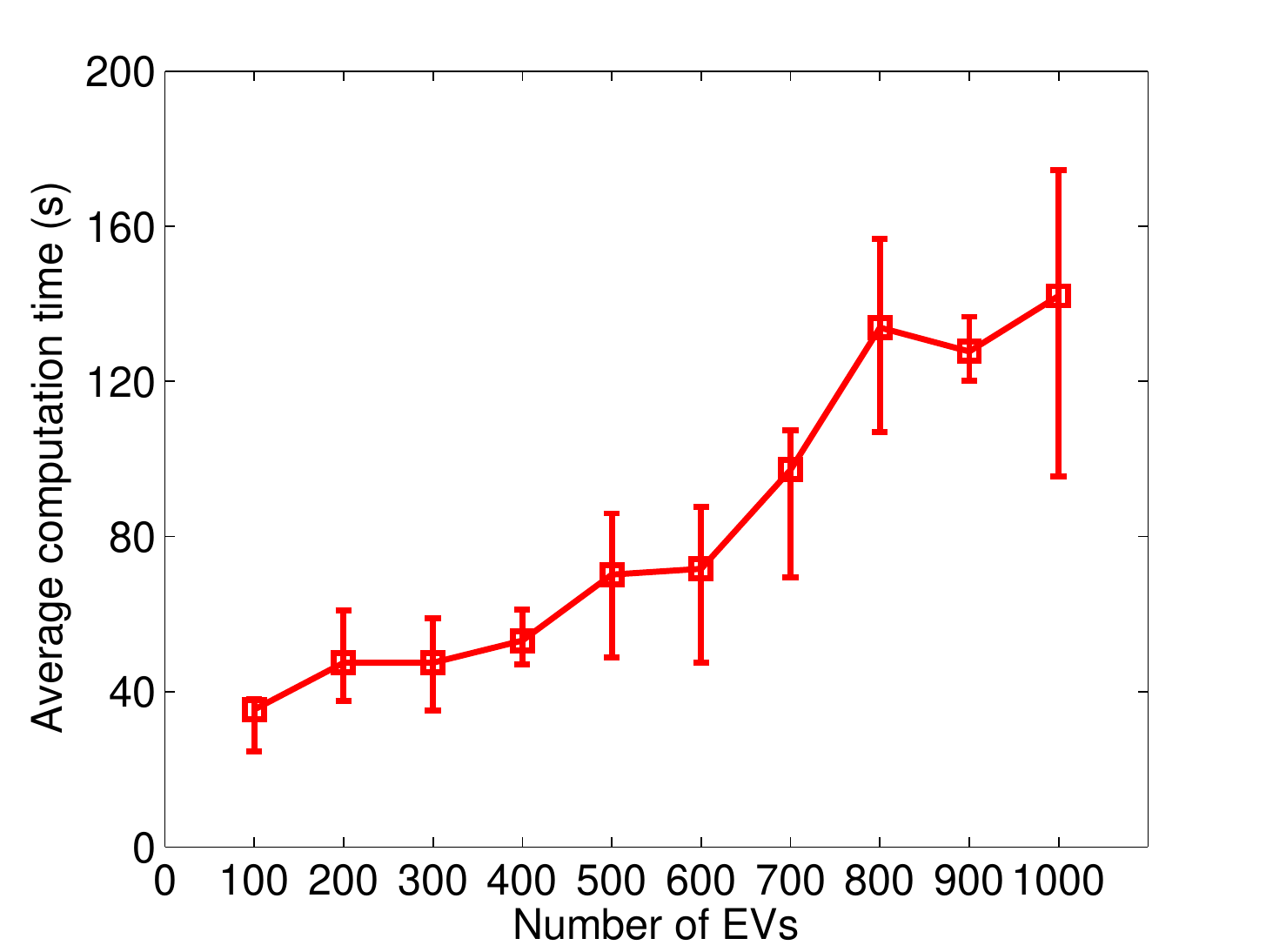}
\caption{Average computation time as a function of \#EVs.}\label{fig:comptime}
\end{figure}

\begin{figure}
\centering
	\includegraphics[width=0.23\textwidth]{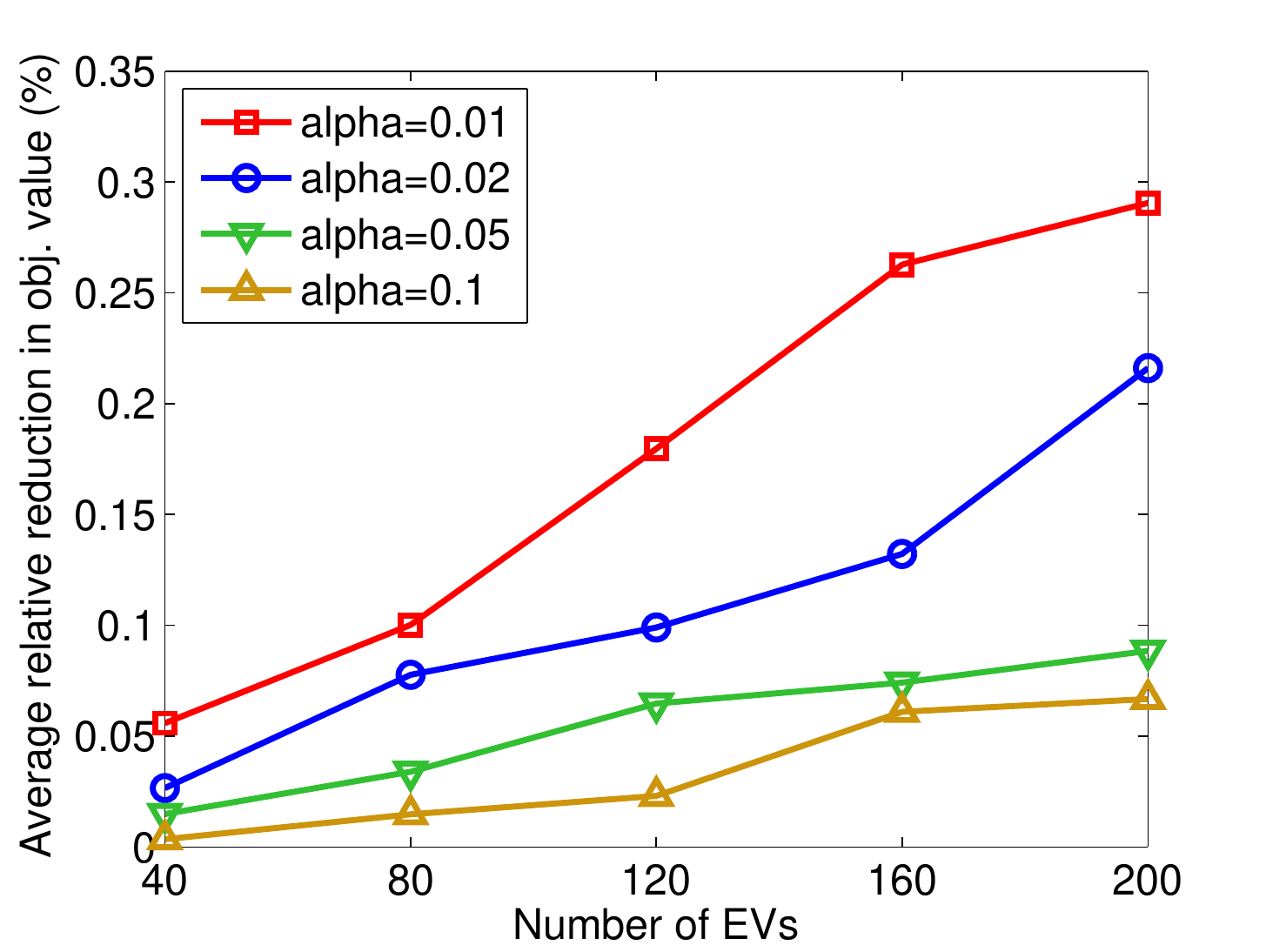}
\caption{Average relative reduction in objective value.}\label{fig:objreduction}
\end{figure}

\begin{figure}
\centering
	\includegraphics[width=0.23\textwidth]{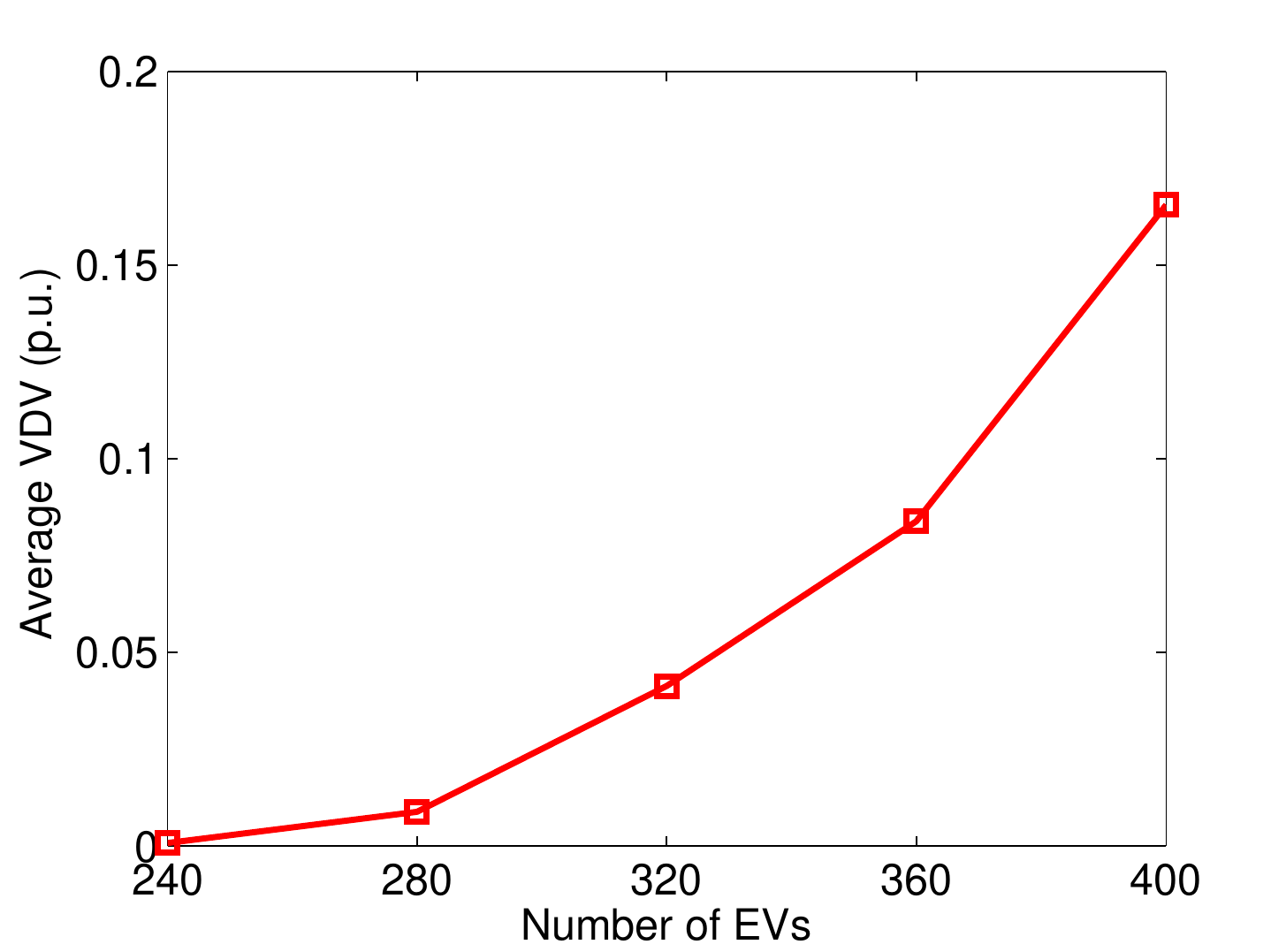}
\caption{Average VDV under nearest-station policy.}\label{fig:VDV}
\end{figure}

\begin{figure}
\centering
\subfigure[]{
	\includegraphics[width=0.21\textwidth]{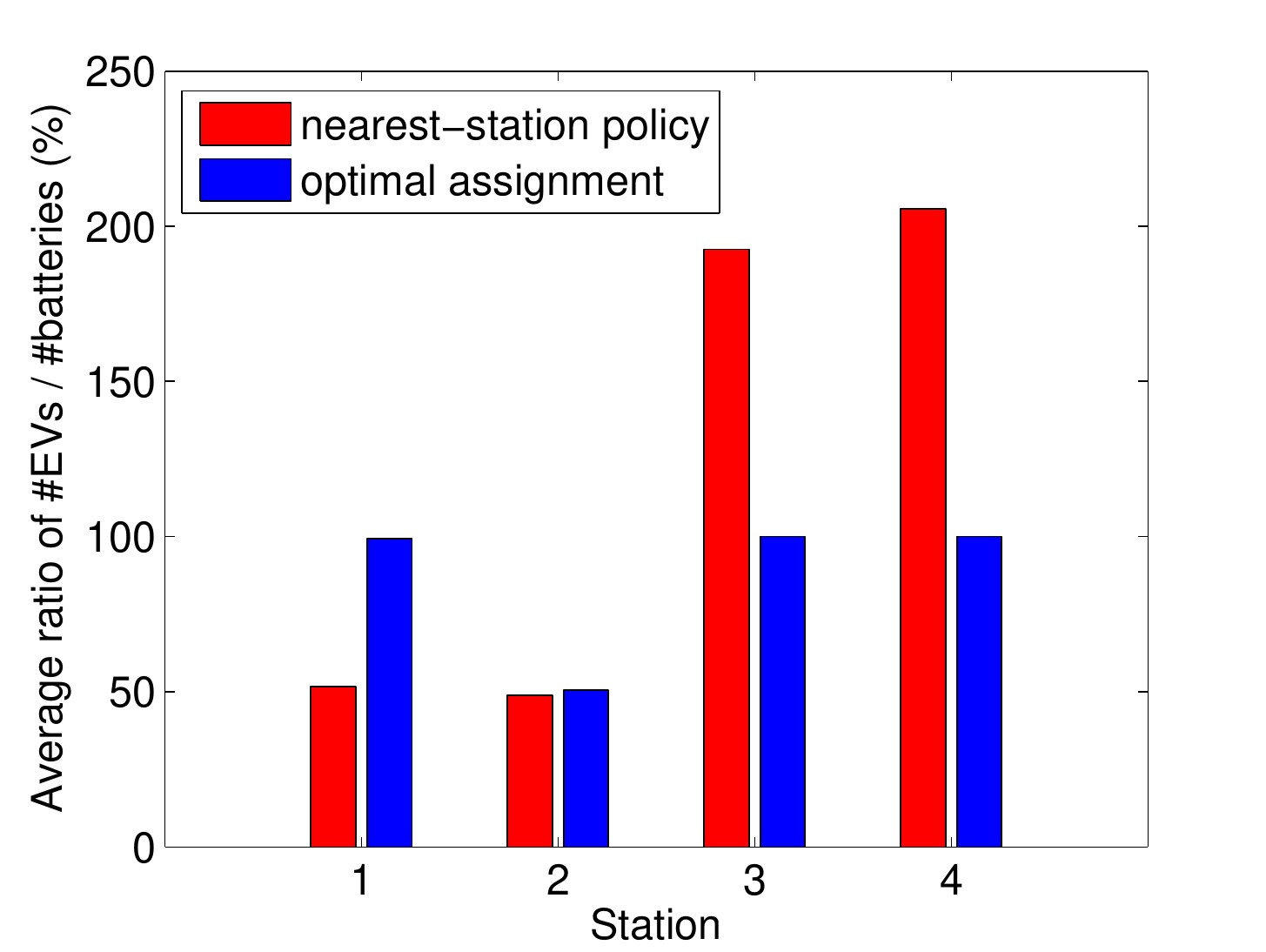}
	\label{fig:prop}
}
\hspace{-0.3in}
\subfigure[]{
	\includegraphics[width=0.21\textwidth]{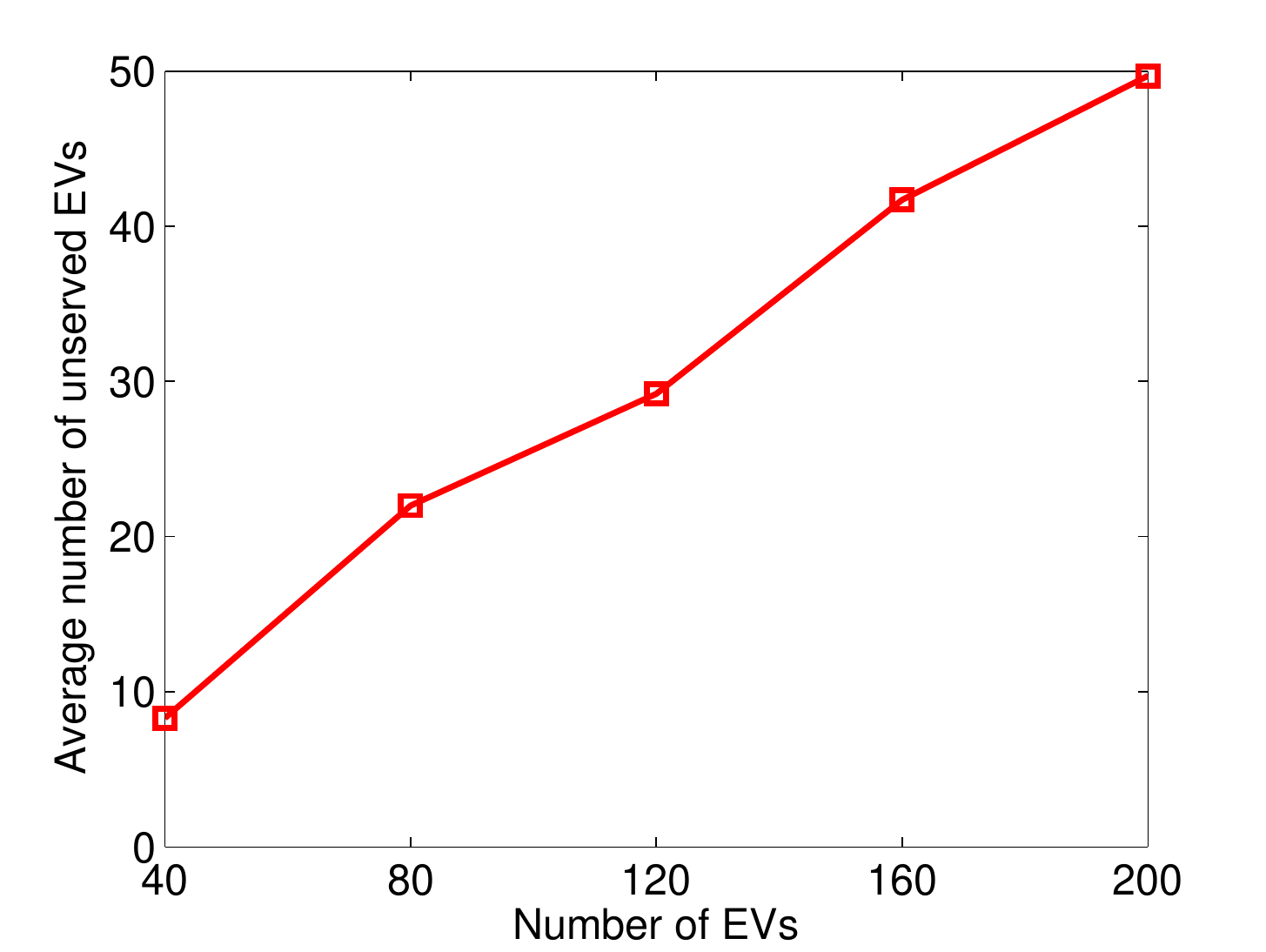}
	\label{fig:unserved}
}
\caption{(a) Average ratio of the number of EVs to that of fully-charged batteries. (b) Average number of unserved EVs under nearest-station policy.}\label{fig:hzissue}
\end{figure}

\noindent
\textbf{Computational effort.}
To demonstrate the potential of the proposed algorithm for practical application, we check its required computational effort by counting its computation time for different number of EVs, since the number of discrete variables in the optimization problem is the computational bottleneck. We use Gurobi to solve the master problem (integer programming) and SDPT3 to solve the slave problem (convex programming) on the MATLAB R2012b platform. Fig.~\ref{fig:comptime}\footnote{Due to the randomization of EVs' initial locations, each datapoint in Fig.~\ref{fig:comptime}, Fig.~\ref{fig:objreduction}, Fig.~\ref{fig:VDV} and Fig.~\ref{fig:hzissue} is an average over 10 simulation runs.} shows the average computation time required by the proposed algorithm to find a global optimum for different numbers of EVs, which validates its computational efficiency.

\noindent
\textbf{Benefit.}
Fig.~\ref{fig:objreduction} displays the average relative reduction in the objective value with different $\alpha$'s using our algorithm, compared with the nearest-station policy.  Scheduling flexibility is enhanced with more EVs, thus improving the savings.  Clearly the smaller the weight $\alpha$ on EVs' travel distance,
the more benefit the proposed algorithm provides over the nearest-station policy.
However, Fig. ~\ref{fig:objreduction} also suggests that the improvement is small, i.e.,  the nearest-station policy is good enough if it is implementable.

The nearest-station policy is sometimes infeasible when there are more EVs nearest to
a station than the number of fully-charged batteries at that station or when some
operational constraints of the distribution network are violated.
In our case study, infeasibility is mainly due to some voltages dropping below their
 lower limits. Define a metric \emph{voltage drop violation} as $\mathrm{VDV}:=\sum_{j\in\mathbb{N}}\max\{\sqrt{\underline{v}_j}-|V_j|,0\}$ to quantify the degree of voltage violation.
Fig.~\ref{fig:VDV} shows the average VDV for the number of EVs ranging from 240 to 400 under the nearest-station policy.
The voltage violation becomes more severe when the number of EVs increases.

It is also interesting to look at cases where there are more EVs nearest to a station than fully-charged batteries that station can provide, which, as far as we know, are common in practice. We reset $M_1=m_1=M_2=m_2=\frac{1}{2}A$ and $M_3=m_3=M_4=m_4=\frac{1}{8}A$ to simulate these situations.
Hence the total number of fully-charged batteries in the system is $\frac{5}{4}A$. Fig.~\ref{fig:prop} shows, for each station, the average ratio of the number of EVs that go to the station for battery swapping to the number of fully-charged batteries at the station, under both the nearest-station policy and an optimal assignment. Under the optimal assignment, 99.40\% of station 1's batteries, 50.60\% of station 2's batteries, and all the batteries at stations 3, 4 are used, thus they have collectively served all $A$ EVs. Under the nearest-station policy, however, only 51.55\% and 48.89\% of stations 1 and 2's batteries respectively (i.e., a total of around $\frac{1}{2}A$ batteries) are used for swapping. At either of stations 3 and 4, the number of EVs is approximately double that of available fully-charged batteries (192.61\% and 205.62\%, respectively). Fig.~\ref{fig:unserved} shows the average number of unserved EVs under the nearest-station policy as a function of the total number $A$ of EVs. On average, a total of $\frac{1}{4}A$ EVs cannot be served at their nearest stations, mainly due to congestion at stations 3 and 4, while available fully-charged batteries at stations 1 and 2 are not fully utilized.

\section{Conclusion}\label{sec:conclusion}

We formulate an optimal scheduling problem for battery swapping that assigns to each EV a best station to swap its depleted battery based on its current location and state of charge. The schedule aims to minimize total travel distance and generation cost over both station assignments and power flow variables, subject to EV range constraints, grid operational constraints and AC power flow equations.
We propose a centralized solution that relaxes the nonconvex constraint of OPF into a second-order cone and then applies generalized Benders decomposition to handle the binary nature of station assignments.
Numerical case studies on the SCE 56-bus distribution feeder show the SOCP relaxation is mostly exact
and generalized Benders decomposition computes an optimal solution efficiently.


\bibliographystyle{ieeetr}
\bibliography{bib}

\begin{thebibliography}{10}

\bibitem{C2ES2016climate}
C2ES, ``{Climate TechBook},'' {\em Center for Climate and Energy Solutions, US:
  www.c2es.org/energy/use/transportation}, 2016.

\bibitem{eia2015annual}
EIA, ``Monthly energy review,'' {\em Energy Information Administration, US
  Department of Energy: www.eia.gov/totalenergy/data/monthly/}, 2015.

\bibitem{leou2014stochastic}
R.-C. Leou, C.-L. Su, and C.-N. Lu, ``Stochastic analyses of electric vehicle
  charging impacts on distribution network,'' {\em IEEE Trans. on Power
  Systems}, vol.~29, no.~3, pp.~1055--1063, 2014.

\bibitem{shang2015orchestrating}
T.~Shang, Y.~Chen, and Y.~Shi, ``{Orchestrating ecosystem co-opetition: Case
  studies on the business models of the EV demonstration programme in China},''
  {\em Electric Vehicle Business Models}, pp.~215--227, 2015.

\bibitem{Low2014a}
S.~H. Low, ``Convex relaxation of optimal power flow, {I}: formulations and
  relaxations,'' {\em IEEE Trans. on Control of Network Systems}, vol.~1,
  no.~1, pp.~15--27, 2014.

\bibitem{Low2014b}
S.~H. Low, ``Convex relaxation of optimal power flow, {II}: exactness,'' {\em
  IEEE Trans. on Control of Network Systems}, vol.~1, no.~2, pp.~177--189,
  2014.

\bibitem{MaCallawayHiskens2013}
Z.~Ma, D.~Callaway, and I.~Hiskens, ``Decentralized charging control of large
  populations of plug-in electric vehicles,'' {\em IEEE Trans. on Control
  Systems Technology}, vol.~21, no.~1, pp.~67--78, 2013.

\bibitem{gan2013optimal}
L.~Gan, U.~Topcu, and S.~H. Low, ``Optimal decentralized protocol for electric
  vehicle charging,'' {\em IEEE Trans. on Power Systems}, vol.~28, no.~2,
  pp.~940--951, 2013.

\bibitem{papadopoulos2013coordination}
P.~Papadopoulos, N.~Jenkins, L.~M. Cipcigan, I.~Grau, and E.~Zabala,
  ``Coordination of the charging of electric vehicles using a multi-agent
  system,'' {\em {IEEE} Trans. on Smart Grid}, vol.~4, no.~4, pp.~1802--1809,
  2013.

\bibitem{han2011estimation}
S.~Han, S.~Han, and K.~Sezaki, ``{Estimation of achievable power capacity from
  plug-in electric vehicles for V2G frequency regulation: Case studies for
  market participation},'' {\em {IEEE} Trans. on Smart Grid}, vol.~2, no.~4,
  pp.~632--641, 2011.

\bibitem{o2014rolling}
A.~O'Connell, D.~Flynn, and A.~Keane, ``Rolling multi-period optimization to
  control electric vehicle charging in distribution networks,'' {\em IEEE
  Trans. on Power Systems}, vol.~29, no.~1, pp.~340--348, 2014.

\bibitem{ChenTong2012}
S.~Chen and L.~Tong, ``{iEMS for large scale charging of electric vehicles:
  Architecture and optimal online scheduling},'' in {\em Proc. of IEEE
  International Conference on Smart Grid Communications (SmartGridComm)},
  pp.~629--634, 2012.

\bibitem{li2014distribution}
R.~Li, Q.~Wu, and S.~S. Oren, ``Distribution locational marginal pricing for
  optimal electric vehicle charging management,'' {\em IEEE Trans. on Power
  Systems}, vol.~29, no.~1, pp.~203--211, 2014.

\bibitem{Tong2014}
Z.~Yu, S.~Li, and L.~Tong, ``On market dynamics of electric vehicle
  diffusion,'' in {\em Proc. of the 52nd Annual Allerton Conferene on
  Communication, Control, and Computing}, pp.~1051--1057, 2014.

\bibitem{sojoudi2011optimal}
S.~Sojoudi and S.~H. Low, ``Optimal charging of plug-in hybrid electric
  vehicles in smart grids,'' in {\em Proc. of IEEE Power \& Energy Society
  General Meeting}, pp.~1--6, 2011.

\bibitem{zhang2015scalable}
L.~Zhang, V.~Kekatos, and G.~B. Giannakis, ``Scalable network-constrained
  electric vehicle charging in multiphase distribution grids,'' {\em arXiv
  preprint arXiv:1510.00403}, 2015.

\bibitem{chen2014electric}
N.~Chen, C.~W. Tan, and T.~Q. Quek, ``Electric vehicle charging in smart grid:
  Optimality and valley-filling algorithms,'' {\em IEEE Journal of Selected
  Topics in Signal Processing}, vol.~8, no.~6, pp.~1073--1083, 2014.

\bibitem{de2015optimal}
J.~de~Hoog, T.~Alpcan, M.~Brazil, D.~A. Thomas, and I.~Mareels, ``Optimal
  charging of electric vehicles taking distribution network constraints into
  account,'' {\em IEEE Trans. on Power Systems}, vol.~30, no.~1, pp.~365--375,
  2015.

\bibitem{7031457}
J.~Franco, M.~Rider, and R.~Romero, ``A mixed-integer linear programming model
  for the electric vehicle charging coordination problem in unbalanced
  electrical distribution systems,'' {\em {IEEE} Trans. on Smart Grid}, vol.~6,
  no.~5, pp.~2200--2210, 2015.

\bibitem{tan2014queueing}
X.~Tan, B.~Sun, and D.~H. Tsang, ``Queueing network models for electric vehicle
  charging station with battery swapping,'' in {\em Proc. of IEEE International
  Conference on Smart Grid Communications (SmartGridComm)}, pp.~1--6, 2014.

\bibitem{yang2014dynamic}
S.~Yang, J.~Yao, T.~Kang, and X.~Zhu, ``{Dynamic operation model of the battery
  swapping station for EV (electric vehicle) in electricity market},'' {\em
  Energy}, vol.~65, pp.~544--549, 2014.

\bibitem{7310884}
P.~You, Z.~Yang, Y.~Zhang, S.~H. Low, and Y.~Sun, ``Optimal charging schedule
  for a battery switching station serving electric buses,'' {\em IEEE Trans. on
  Power Systems}, vol.~31, no.~5, pp.~3473--3483, 2016.

\bibitem{sarker2015optimal}
M.~R. Sarker, H.~Pand{\v{z}}i{\'c}, and M.~A. Ortega-Vazquez, ``Optimal
  operation and services scheduling for an electric vehicle battery swapping
  station,'' {\em IEEE Trans. on Power Systems}, vol.~30, no.~2, pp.~901--910,
  2015.

\bibitem{zheng2014electric}
Y.~Zheng, Z.~Y. Dong, Y.~Xu, K.~Meng, J.~H. Zhao, and J.~Qiu, ``Electric
  vehicle battery charging/swap stations in distribution systems: comparison
  study and optimal planning,'' {\em IEEE Trans. on Power Systems}, vol.~29,
  no.~1, pp.~221--229, 2014.

\bibitem{zhang2016benefit}
X.~Zhang and R.~Rao, ``A benefit analysis of electric vehicle battery swapping
  and leasing modes in china,'' {\em Emerging Markets Finance and Trade},
  vol.~52, no.~6, pp.~1414--1426, 2016.

\bibitem{substationserviceradius}
Q.~Wang, K.~Qin, and D.~Chen, ``Power supply radius optimization for
  substations with large capacity transformers.''
  \url{http://www.cqvip.com/read/read.aspx?id=83887068504849524852484852}.

\bibitem{baran1989optimal}
M.~E. Baran and F.~F. Wu, ``Optimal capacitor placement on radial distribution
  systems,'' {\em IEEE Trans. on Power Delivery}, vol.~4, no.~1, pp.~725--734,
  1989.

\bibitem{benders1962partitioning}
J.~F. Benders, ``Partitioning procedures for solving mixed-variables
  programming problems,'' {\em Numerische mathematik}, vol.~4, no.~1,
  pp.~238--252, 1962.

\bibitem{geoffrion1972generalized}
A.~M. Geoffrion, ``Generalized benders decomposition,'' {\em Journal of
  optimization theory and applications}, vol.~10, no.~4, pp.~237--260, 1972.

\bibitem{7287792}
P.~You, Z.~Yang, M.~Y. Chow, and Y.~Sun, ``Optimal cooperative charging
  strategy for a smart charging station of electric vehicles,'' {\em IEEE
  Trans. on Power Systems}, vol.~31, no.~4, pp.~2946--2956, 2016.

\bibitem{farivar2012optimal}
M.~Farivar, R.~Neal, C.~Clarke, and S.~Low, ``{Optimal inverter var control in
  distribution systems with high PV penetration},'' in {\em Proc. of IEEE Power
  \& Energy Society General Meeting}, pp.~1--7, 2012.

\bibitem{yilmaz2013review}
M.~Yilmaz and P.~T. Krein, ``Review of battery charger topologies, charging
  power levels, and infrastructure for plug-in electric and hybrid vehicles,''
  {\em IEEE Trans. on Power Electronics}, vol.~28, no.~5, pp.~2151--2169, 2013.

\bibitem{eia2011annual}
EIA, ``Annual energy review,'' {\em Energy Information Administration, US
  Department of Energy: www.eia.doe.gov/emeu/aer}, 2011.

\end{thebibliography}

\end{document}